\documentclass[11pt]{article}
\usepackage{amssymb}
\usepackage{amsthm}
\usepackage{amsfonts}
\usepackage{mathtools}
\usepackage{amsmath}
\usepackage{bbm}
\usepackage{multirow}
\usepackage{url}
\usepackage{breakurl}
\usepackage[hidelinks]{hyperref}

\usepackage{algorithmic}
\usepackage{algorithm}
\usepackage{setspace}
\usepackage[font=footnotesize]{subcaption}
\usepackage[font={footnotesize, stretch=1.5}]{caption}
\usepackage[normalem]{ulem}
\usepackage{soul}
\usepackage{graphicx}
\usepackage[usenames,dvipsnames]{xcolor}
\usepackage{epstopdf}

\newtheorem{theorem}{Theorem}[section]

\newcommand*{\email}[1]{\href{mailto:#1}{\nolinkurl{#1}} }

\title{Optimal policy for control of epidemics with constrained time intervals and region-based interactions}
 \author{Xia Li\thanks{ Microsoft, WA (Corresponding author: \email{xli51@g.ucla.edu} or \email{xiali1@microsoft.com})}
 \and Andrea L. Bertozzi \thanks{Departments of Mathematics and Mechanical and Aerospace Engineering, University of California, Los Angeles (\email{bertozzi@math.ucla.edu})}
 \and  P. Jeffrey Brantingham \thanks{Department of Anthropolog, University of California, Los Angeles (\email{branting@g.ucla.edu})}
 \and Yevgeniy Vorobeychik \thanks{Department of Mathematics, Washington University in St. Louis (\email{yvorobeychik@wustl.edu})}
 \thanks{The code for all the experiments and the data can be found: \url{https://github.com/xli07/PolicySIR/tree/main}}
 }
 
\begin{document}



\date{}
\maketitle

\begin{abstract}
We introduce a policy model coupled with the 
susceptible–infected–
recovered (SIR) epidemic model to study interactions between policy-making and the dynamics of epidemics. 
We consider both single-region policies, as well as game-theoretic models involving interactions among several regions, and hierarchical interactions among policy-makers modeled as multi-layer games.
We assume that the policy functions are piece-wise constant with a minimum time interval for each policy stage, considering policies cannot change frequently in time or they cannot be easily followed. The optimal policy is obtained by minimizing a cost function which consists of an implementation cost, an impact cost, and, in the case of multi-layer games, a non-compliance cost. We show in a case study of COVID-19 in France that when the cost function is reduced to the impact cost and is parameterized as the final epidemic size, the solution approximates that of the optimal control in  Bliman {\em et al}, J. Optim. Theory Appl., 189, 2021, for sufficiently small minimum policy time interval.
For a larger time interval however the optimal policy is a step down function, quite different from the step up structure typically deployed during the COVID-19 pandemic.  In addition, we present a counterfactual study of how the pandemic would have evolved if herd immunity was reached during the second wave in the county of Los Angeles, California. Lastly, we study a case of three interacting counties with and without a governing state.  
\end{abstract}


\section{Introduction}	
In the course of battling COVID-19, public health policies sought to enforce non-pharmaceutical interventions to slow or halt the spread of the pandemic. Common policies included `safer-at-home', `social distancing' and `mask wearing' mandates, which were seen as crucial during the early stages of the pandemic prior to the availability of vaccines. The timeline of COVID-19 globally and locally (\cite{CDC_timeline, ca_timeline}) indicates that the evolution of policy affected the evolution of the pandemic and vice versa. For example, in the county of Los Angeles, social distancing was first mandated \cite{department_of_public_health} on March 21, 2020, about a month after the first reported COVID-19 case in LA. Around that time, the Los Angeles Mayor's Office released the `safer-at-home' policy \cite{office_city_of_los_angeles}. One week later, beaches, hiking trails, dog parks, skate parks, and other public sites and facilities were temporarily closed. On April 15th, as infected cases continued to increase, facial coverings were mandated in many indoor places \cite{lin_2020}. In hindsight, it is important to ask: Were policies that were enforced done so in an optimal way? What can we learn by using mathematical modeling to understand the interplay between policy and spread of disease?  This paper introduces a policy model coupled to a susceptible–infected–recovered (SIR) epidemic model to study interactions between policy-making and the dynamics of epidemics. 
There have been several studies on the relationship between  policies and epidemics \cite{blower2002health,pnas_1918, Allison2018, burke2006,ijerph19106119}. In a study analyzing data from 16 US cities during the 1918 pandemic \cite{pnas_1918}, Bootsma and Ferguson analyzed specific outcomes related to the impact of the delay of lockdown policies on the total deaths and also on the appearance of second waves of outbreaks due to reopening too early. The analysis was done fitting available data to an SEIR model. They also considered optimal control for the simpler SIR and the end-state of the pandemic, noting that there exists an optimal control level with fewer deaths and no second wave. More recently, Bliman et al. \cite{bliman_2021_optimal} developed a theoretical study of the optimal control of a classical SIR outbreak.  
Bliman et al. do not consider the possible of vaccines or pharmaceutical interventions.  Rather, focusing exclusively on non-pharmaceutical interventions, they design an optimal policy that achieves an end state as close to the herd immunity threshold as possible. This is the same problem considered briefly in a section of \cite{pnas_1918}.  Bliman et al. prove the existence and the uniqueness of their solution and showed the optimal social distancing polity is a bang-bang controller \cite{sd_reduce_fz},  generalizing the results of  \cite{bliman_2021_optimal} by modeling without prescribing the starting date of the policy.

The substantial theoretical insights of Bliman et al.'s model are limited in their practical implications by a few assumptions. First, Bliman et al. assume that policy that can change continuously in time, which would imply, for example, the ability to shift in three successive instants between no restrictions, perfect ``lockdown", and back to no restrictions. 
As observed during the COVID-19 pandemic, policies that change frequently in time cannot be easily followed. Moreover, policies must be relatively easy to interpret, with a small number of different intensity levels (see Fig. \ref{fig:ca_p_real}). A practical implementation also requires a minimum time duration for a particular stage of the policy. These practical constraints can be modeled together as a piece-wise constant function of time with a minimum time interval for each well-defined policy level (i.e., not continuous).  With this idea in mind, we aim to re-examine the optimal practical policy among all possible piece-wise constant policies with minimal time duration.

Second, Bliman et al. assume that that  the only outcome to manage is the final epidemic size. This so-called ``impact cost" is clearly a central concern (see below). However, as also seen during the pandemic, there are real trade-offs between decreased infections and the negative impact of strict policies on other aspects of society such as remote learning for young students, employment curtailment in certain job sectors, and lack of key services provided to the public. In the present work, we modify Bliman et al.'s model to take into account these other practical ``implementation costs." Specifying a short minimal time interval during which policies must remain constant (e.g., one week), we find our results resemble Bliman et al.'s bang-bang controller \cite{bliman_2021_optimal} despite the more complex cost structure that includes both impact and implementation costs. With a larger minimal time interval during which time policies must remain constant (e.g., 28 days), optimal policies depart from the bang-bang solution.

Finally, Bliman et al. also assume a pandemic spreading in an single population pool overseen by a single policy-making entity. The reality of the COVID-19 pandemic is that there are policy makers at several (nested) hierarchical scales that oversee different population pools. For example, within the United States, policies may be set at Federal, State, County and local levels, not to mention finer-grained institutional and family scales. And populations at any one scale (e.g., counties) may interact to varying degrees. Inspired by the work of Jia et al. \cite{jia2021game}, we introduce a hierarchical version of Bliman et al.'s model with sequential (Stackleberg) policy-making. Specifically, levels higher in a jurisdictional hierarchy make policy decisions, while levels lower in the hierarchy make their decisions with full knowledge of the policy recommendations from above. We find that a hierarchical structure can make the policies converge in all regions using the right weight for a non-compliance cost. 

The remainder of this paper is organized as follows. We first introduce the work in \cite{bliman_2021_optimal} and reproduce the results using our methods. We discuss how different optimal policies result from different parameter choices for model constraints and costs. Next, we discuss an empirical case study of the so-called  ``second wave" of the pandemic (November 6th 2020–May 12th 2021) in Los Angeles County, California. Last, we use simulation to study optimal control of the pandemic in three counties with and without a governing state as an example of the multi-layer multi-regional model.

\section{Policy model using optimal control}
A policy function is a continuous function that has a range of $[0,1]$. As the numerical value increases, the strictness of the policy decreases. The Numerical value $0$ denotes a total lockdown and $1$ denotes no control. We assume a policy $u(t)$ directly influences the level of a lockdown, which affects the rate of the population transport from compartment $S$ to $I$. We use the following policy-incorporated SIR:
\begin{eqnarray}\displaystyle\label{eqn:sir_ode}
\begin{cases}
&\dfrac{d S(t)}{dt} =  - u(t) \beta \frac{I(t)S(t)}{N} , \\
&\dfrac{dI(t)}{dt} = u(t)\beta \frac{I(t)S(t)}{N} -\gamma  I(t) , \\
&\dfrac{dR(t)}{dt} =  \gamma  I(t) , \\
&  S(0)=    S_0,  \quad  I(0)=  I_0, \quad   R(0)=  R_0.
\end{cases}
\end{eqnarray}
Like the traditional SIR model, the reproduction number $\mathcal{R}_0 = \frac{\beta}{\gamma}$. Herd immunity occurs when a large proportion of the population has become immune to the infection. Mathematically, it is defined as the value of $S$ below which the number of infected decrease and can be calculated as $S_{\text{herd}} = \frac{N}{\mathcal{R}_0}$. In \cite{bliman_2021_optimal}, 
a policy $u(t)$ is assumed to belong to the admissible set $\mathcal{U}_{\alpha_{\max},T_0}$ defined by
$$\{ u\in \mathcal{L}^{\infty}([0,+\infty]), \alpha_{\max} \leq u(t) \leq 1 \text{ if } t \in[0,T_0], u(t) = 1 \text{ if } t > T_0\}.$$
The constant $T_0$ characterizes the duration of the policy, and $\alpha_{\max}$ its maximal intensity. 
In \cite{bliman_2021_optimal}, Theorem 2.1 states that no finite time intervention is able to stop the epidemics before or exactly at the herd immunity. However, one may stop arbitrarily close to herd immunity by having a sufficiently long intervention of sufficient intensity. To determine the closest state $S$ to this threshold attainable by control of maximal intensity $\alpha_{\max}$ on the interval $[0, T_0]$, one is led to consider the following optimal control problem:
\begin{equation}\label{eqn:oc}
    \sup_{u\in\mathcal{U}_{\alpha_{\max},T_{0}}}S_{\infty}(u).
\end{equation}
Furthermore, Bliman et al. prove the existence and uniqueness of the optimal solution to problem \ref{eqn:oc} and that the solution is a bang-bang controller (a control that switches from one extreme to the other). More specifically, they have the following theorem:
\begin{theorem} (Theorem 2.1 in \cite{bliman_2021_optimal})
Let $\alpha_{\max} \in [0, 1)$ and $T_0 > 0$. Problem \ref{eqn:oc} admits a unique solution $u^*$.
Furthermore,\\
(i) the maximal value $S^*_{\infty,\alpha_{\max},T_0} \coloneqq \{ \max{S_{\infty}(u) : u \in \mathcal{U}_{\alpha_{\max},T_0}}\} $is non-increasing with
respect to $\alpha_{\max}$ and non-decreasing with respect to $T_0$ .\\
(ii) there exists a unique $T_0^* \in [0, T_0 )$ such that $u^* = u_{T_0^*}\coloneqq \mathbbm{1}_{[0,T_0^*]}+\alpha_{\max}\mathbbm{1}_{[T_0^*,T_0]}+\mathbbm{1}_{[T_0,+\infty)}$ (in particular, the optimal control is bang-bang).
\end{theorem}

\section{Single region case}\label{sec:single_region}
We use the same policy-incorporated SIR model for the epidemic dynamic as in \cite{bliman_2021_optimal}. Instead of minimizing the final epidemic size alone, we adopt a similar policy-making process as in \cite{jia2021game} by using a cost function that takes into account the cost of implementing the policy, the impact of the infection and a penalty for being non-compliant. The latter cost only applies in hierarchical models where a lower-level unit can choose to not follow the policy recommendation of a higher-level unit.

\begin{figure}[!h]
    \centering
    \includegraphics[width = 0.99\linewidth]{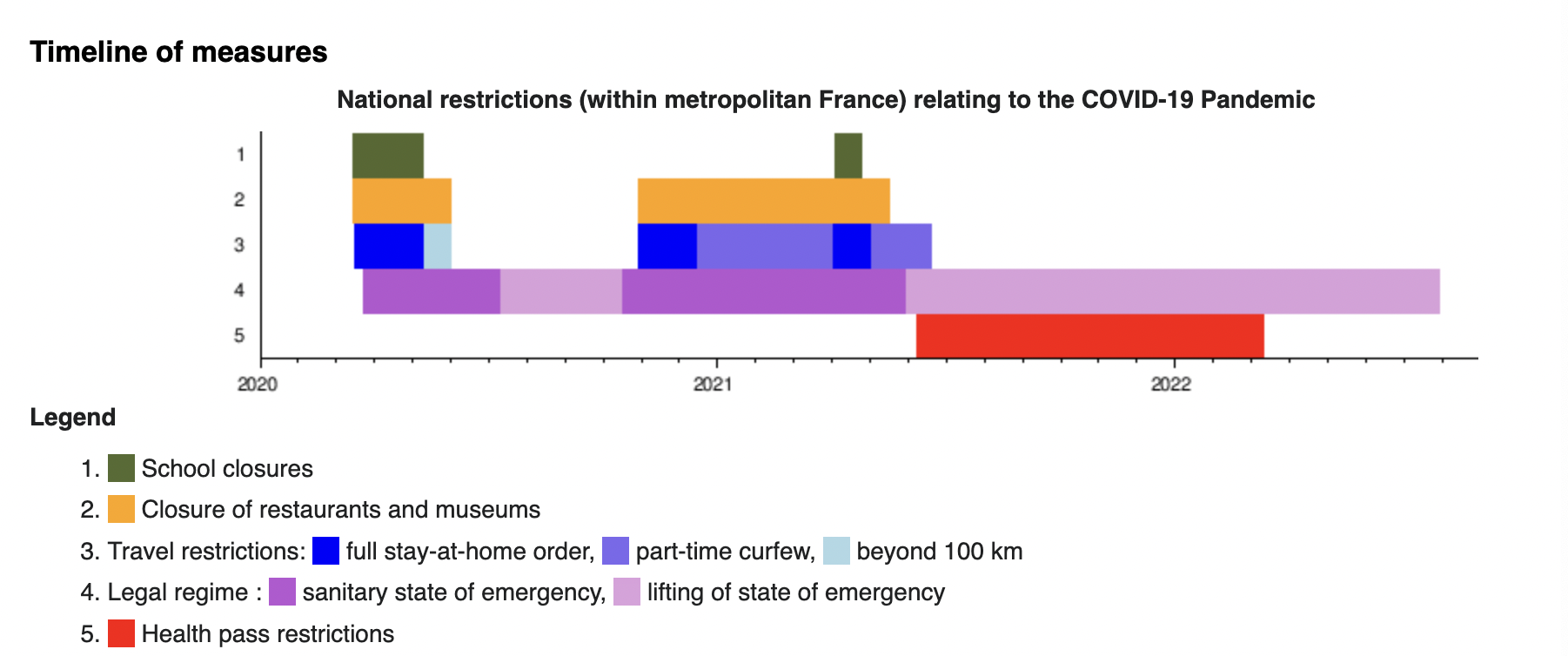}
    \caption{Timeline of COVID-19 restrictions in France.}
    \label{fig:France-timeline}
\end{figure}

We also consider practical implementation constraints, namely that the policy can only be implemented using a finite number of discrete levels of control and with a minimal time interval during which a policy must remain constant. As an example, consider the policy implementation in France during the year 2020 and 2021 shown in Fig.~\ref{fig:France-timeline} (\cite{wikipedia_2022}). Implemented policies were discrete both in terms of the small number of intervention types and the fixed time intervals of enforcement, the shortest of which was approximately 15 days in duration, with the longest lasting more than a year. A discrete policy model is realistic given the empirical pattern of real-world interventions. Such a model also simplifies the computation problem of optimal policy discovery by searching through a discrete set of potential policies rather than a continuum of policies.
\subsection{The policy-incorporated SIR model}
To model the evolution of the pandemic, we discretize the system of ODE using forward Euler's method with a time step of 1:
\begin{eqnarray}\displaystyle
\begin{cases}
\label{eqn:dis}
& S(t) = S(t-1) - \alpha \beta \frac{I(t-1)S(t-1)}{N} , \\
&I(t) = I(t-1) + \alpha\beta \frac{I(t-1)S(t-1)}{N} -\gamma  I(t-1), \\
&R(t)= R(t-1) + \gamma  I(t-1) , \\
&  S(0)=    S_0,  \quad  I(0)=  I_0, \quad   R(0)=  R_0.
\end{cases}
\end{eqnarray}
Where $\alpha = u(t-1)$. Equation \eqref{eqn:dis} can be seen as a first order approximation of the system of ODE in \eqref{eqn:sir_ode}.  \\
\subsection{The policy model}
\paragraph{Policy function}\label{sec:policy_fun}
Instead of assuming continuous policy functions, we consider a more realistic set of policies with a discrete number of different stages and intensity levels. Therefore, policy functions form a subset of the admissible set $\mathcal{U}_{\alpha_{\max},T_0}$ in \cite{bliman_2021_optimal}.

We define the minimal policy time interval (MPTI) as the minimal duration time during which a policy remains constant or unchanged. This notion assumes that there is a minimal duration time of different stages of a policy. In addition to $u \in \mathcal{U}_{\alpha_{\max},T_0}$, we assume that every policy $u$ has a minimal policy time interval $\Delta t$ and in our simulations, the duration of each stage is a multiple of the MPTI. We denote this subset of policy functions as $\mathcal{U}_{\alpha_{\max},T_0}^{\Delta t}$. In the past, many public health agencies enforced policies for time periods that corresponded to the work week (e.g. seven days) or multiples of this (e.g. one month). For the purpose of this paper, we assume the MPTI is an integer multiple of seven days.

We additionally assume that policy functions take values from a finite number of intensity levels $A \subset[\alpha_{\max},1]$, corresponding to different stages of the policy. In the simulations, we use $ A = \{\alpha_{\max}, \frac{\alpha_{\max}+1}{2},1\}$. As a result, the policies we consider are piece-wise functions.

\paragraph{Cost function} 
At time $t$, let $u(t) = \alpha$. The cost at time $t$ is defined by:
\begin{equation}\label{eqn:cost}
    c(\alpha) = \kappa c^{\text{implementation}}(\alpha) + \eta c^{\text{impact}}(\alpha) + (1-\kappa-\eta)c^{\text{non-compliance}} (\alpha)
\end{equation}
The cost function is a linear combination of three parts: 
\begin{enumerate}
    \item the implementation cost, which represents the consequences of policies meant to curtail the pandemic on individuals and the broader economic and social systems.
    \item the impact cost, which represents the consequences of people getting sick both on individuals and the broader economic and social systems.
    \item the non-compliance cost, which is a penalty imposed by a policy-maker upon an agent within its jurisdiction for deviating from its recommendation (e.g., a fine or litigation costs).
\end{enumerate}
The implementation cost is a non-increasing function of $\alpha$ and the impact cost function is a non-decreasing function of $\alpha$.  The coefficients $\kappa,\eta,\kappa + \eta \in [0,1]$. 
The cost from time $t_1$ to $t_2$ is defined as the averaged integral of the cost function over a total time period $T$:
\begin{equation}
    c_{t_1t_2}(u)
    =\frac{1}{T}\int_{t_1}^{t_2}c(\alpha(t)) dt .
\end{equation}
There are different ways to parameterize the cost function. In this paper, the cost function is parameterized in the following way:
\begin{equation}\label{eqn:para_cost}
    c_{t_1t_2}(u)
    = \kappa\left(1-\frac{\int_{t_1}^{t_2 }u(t) dt}{T}\right) + \eta R_{t_2}(u) + (1-\kappa-\eta)\frac{1}{T}\int_{t_1}^{t_2}(u(t) - \pi (u(t)))^2 dt ,
\end{equation}
where $R_{t_2}(u)$ is the fraction of the recovered population at time $t_2$ if policy $u$ is adopted during $[t_1,t_2]$ and $\pi(u)$ is the policy of the agent one level above. The parameterization of the implementation cost and the non-compliance cost are adopted from \cite{jia2021game}. The impact cost is parameterized as the recovered population at time $t_2$ to approximate the impact on the medical system since a fraction of the recovered represents the hospitalized population. If the cost function $u$ is fixed at  constant value $\alpha$ over time interval $[t_1,t_2]$ , the cost can be written as:
\begin{equation}
    c_{t_1t_2}(u)
    = \kappa\frac{t_2-t_1}{T}(1-\alpha) + \eta R_{t_2}(\alpha) + (1-\kappa-\eta)\frac{t_2-t_1}{T}(\alpha - \pi (\alpha))^2.
\end{equation}
An example of cost functions with different weights using the above parameterization is shown in Fig. \ref{fig:cost}. 
In our simulation for a single region, we use a averaged total cost over a time period $T$ as the following:
\begin{equation}\label{eqn:total_cost}
\begin{aligned}
    c_{\text{total}} (u)
    &=\frac{1}{T}\int_{0}^{T}c(u(t)) dt \\
    &= \kappa\left(1-\frac{\int_{0}^{T}u(t) dt}{T}\right) + \eta R_{T}(\alpha) + (1-\kappa-\eta)\int_{0}^{T}(u(t) - \pi (u(t)))^2 dt. 
\end{aligned}
\end{equation}
If at time $T$, the SIR model has reached the equilibrium, we can use $ R_{T}(\alpha)$ to approximate $R_{\infty}$, the fraction of the final size of the recovered population. To find the optimal policy, we solve for the following optimization problem:
\begin{equation}\label{eqn:prob}
 u(t) = \arg\min_{u'}   c_{\text{total}} (u') 
\end{equation}

\subsection{Algorithm}
We discretize time by MPTI $\Delta t$ and the policy intensity into multiple levels. Let $T$ be the total time and $A$ be the set of possible policy intensities (e.g., $A=\{0,0.5,1\}$). 
We search for all the policies that lead to $S_\text{final}$ being close to $S_\text{herd}$, i.e. $ S_\text{final} > S_\text{herd} - \epsilon$, for some sufficiently small $\epsilon$ using a depth-first search algorithm \cite{dfs}. The depth-first search algorithm stores the cost up to the current time interval and reuses this result to obtain the total cost for each policy function through backtracking. Let $N = \frac{T}{\Delta t}$ and $N$ denote the number of stages of a policy. In total, there are $|A|^N$ policies. We initialize the minimal cost $c_{\min}$ to be 9999. Assume the initial susceptible and infected population are $S_0$, $I_0$, respectively. For $n$-th time interval ($n < N$), we choose a value from the set intensity levels $A$ that has not been used before, calculate the cost for the policy intensity, add it to the previous cost, and calculate the susceptible and the infected at the end of $n$-th time interval using the chosen intensity. Then we move to $(n+1)$-th time interval. If the end time interval is reached, we check if  $S_\text{final} > S_\text{herd} - \epsilon$. If yes, we calculate the cost for the final time interval and add it to the previous cost to get the current total cost $c$. If the total cost $c$ is smaller than $c_{\min}$, we update $c_{\min}$ with the total cost $c$, and the optimal policy $u_\text{opt}$ with $u$. Next, we go back to the previous time interval and repeat the same procedure.  After searching over all policies, the policy with the lowest cost is the optimal policy. The detailed algorithm is presented in Alg.~\ref{alg:single_psir}. 
\begin{algorithm}[!h]
\caption{\textsc{Single-region policy SIR}}\label{alg:single_psir}
    \begin{algorithmic}[1]
        \STATE \textbf{Input}: \textbf{ Time} $T$,\textbf{ initial infected population} $I_0$, \textbf{ initial susceptible population} $S_0$, \textbf{ intensity levels} $A$, \textbf{ minimal policy time interval} $\Delta t$,\textbf{ policy end time} $T_0$, \textbf{ Tol} $\epsilon$
        \STATE Initialize \textbf{ county policies},\textbf{ minimal cost $c_{\min} = 9999$}, \textbf{ current cost $c = 0$}
        \STATE $N = \frac{T}{\Delta t}$, $n=1$
        \STATE Initialize \textbf{ policy} $u\in\mathbb{R}^{N}$, \textbf{ optimal policy} $u_{\text{opt}}\in\mathbb{R}^{N}$
        \IF{$n==N$}\label{line:4}
            \FOR{intensity level $\alpha \in A$}
                \STATE calculate $S_{\text{final}}$ using the intensity level $\alpha$, $S_{N-1}$, $I_{N-1}$, and the update rule (\ref{eqn:dis})
                \IF{$S_\text{final} > S_\text{herd} - \epsilon$}
                    \STATE calculate the cost $c_{\text{temp}}=C(\alpha)$ for $N$-th time interval
                    \STATE $c \mathrel{+}= c_{\text{temp}}$,$u(N) = \alpha$
                     \IF {$c\leq c_{\min}$}
                            \STATE $c_{\min} = c$, $u_{\text{opt}} = u$ 
                    \ENDIF
                    \STATE $c \mathrel{-}= c_{\text{temp}}$
                \ENDIF
            \ENDFOR
        \ELSE
            \FOR{intensity level $\alpha \in A$}
                \STATE calculate the cost $c_{\text{temp}}=C(\alpha)$ for the $n$-th time interval
                \STATE calculate the susceptible $S_n$ and the infected $I_n$ at time $n\Delta t$ using $\alpha$, $S_{n-1}$, $I_{n-1}$ and the update rule \refeq{eqn:dis}
                \STATE $c \mathrel{+}= c_{\text{temp}}$, $u(n) = \alpha$\footnotemark
                \STATE $n+=1$ \label{line:19}
                \STATE repeat line \ref{line:4}–\ref{line:19}  until $n=N$
                \STATE $n-=1$
                \STATE $c \mathrel{-}= c_{\text{temp}}$
            \ENDFOR 
            
        \ENDIF
    \RETURN $c_{\min},u_{\text{opt}}$
    \end{algorithmic}
\end{algorithm}
\footnotetext{$u(i)$ represents the $i$-th entry of vector $u$.}
\subsection{Simulations}
In this section, we present the results for both single-region and multiple-region cases. We first compare the results of our discretized method of the COVID-19 in France with the results  \cite{bliman_2021_optimal}. Next, we study the second wave (November 2020–May 2021) in Los Angeles County. 

\subsubsection{Optimal policy in France}
\begin{figure}[!ht]
    \centering
    \begin{subfigure}[t]{0.45\linewidth}
        \centering
         \includegraphics[width=\linewidth]{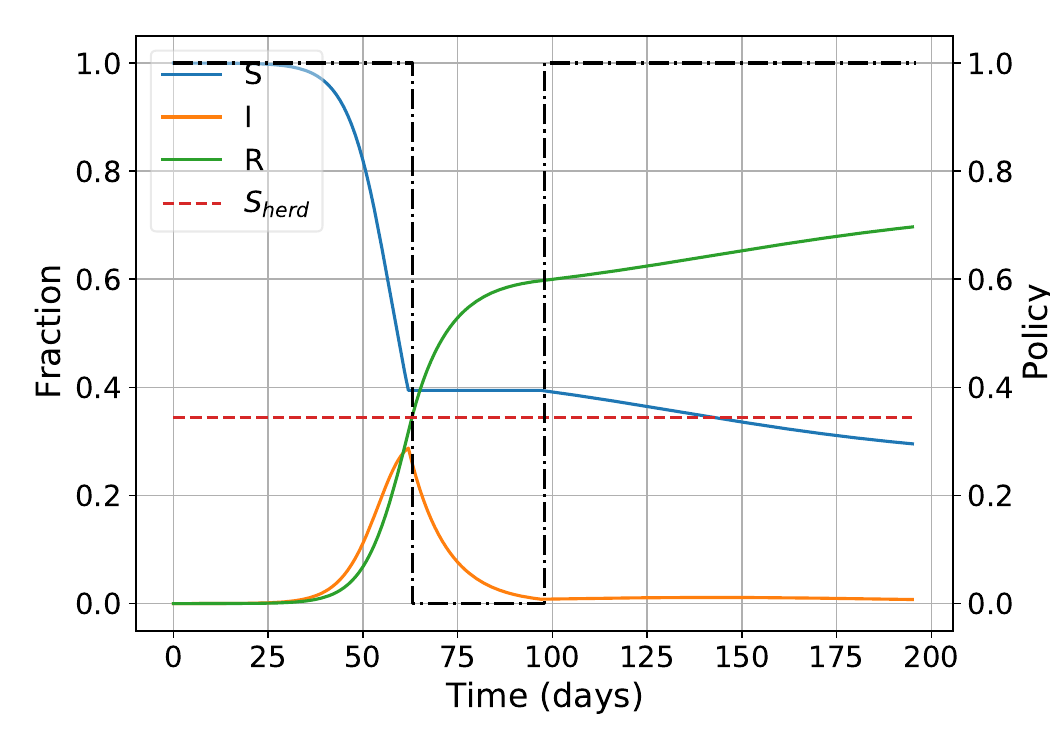}
         \caption{Optimal policy with the minimal policy time interval $\Delta t= 7$ days, $S_{\infty} = 0.296$}
         \label{fig:singcounty_dt=5}
    \end{subfigure}
    \hspace{1em}
    \centering
    \begin{subfigure}[t]{0.45\linewidth}
        \centering
         \includegraphics[width=\linewidth]{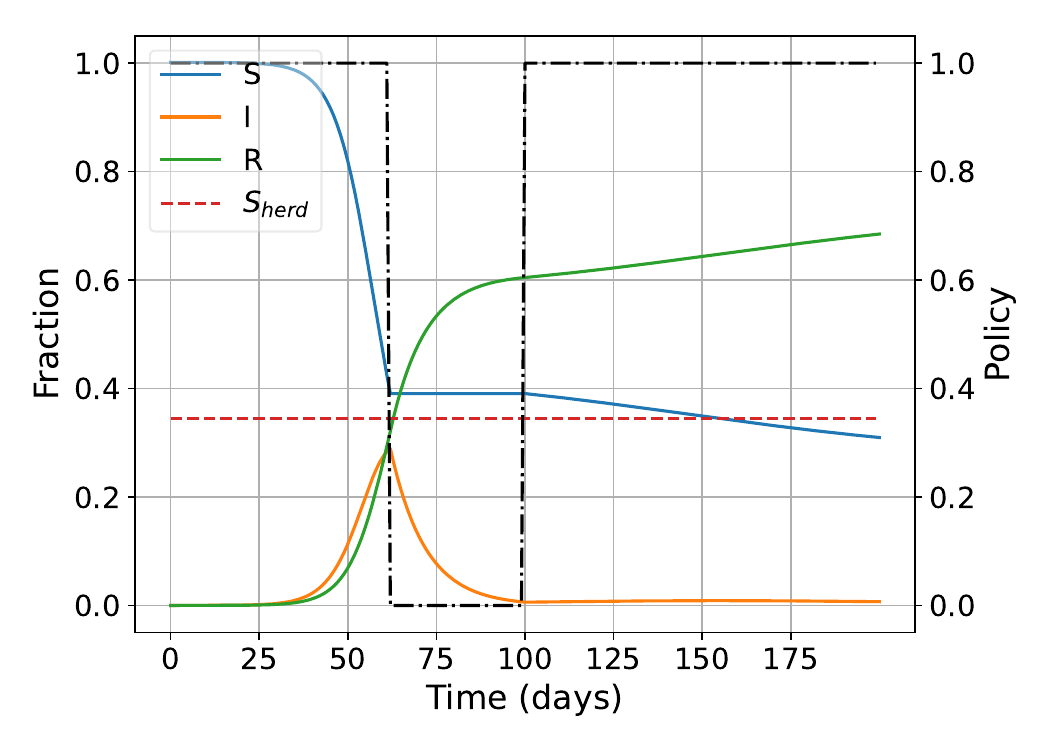}
         \caption{Optimal policy in \cite{bliman_2021_optimal}, $S_{\infty} = 0.31$}
         \label{fig:bliman}
    \end{subfigure}
    \\
    \centering
    \begin{subfigure}[t]{0.45\linewidth}
        \centering
         \includegraphics[width=\linewidth]{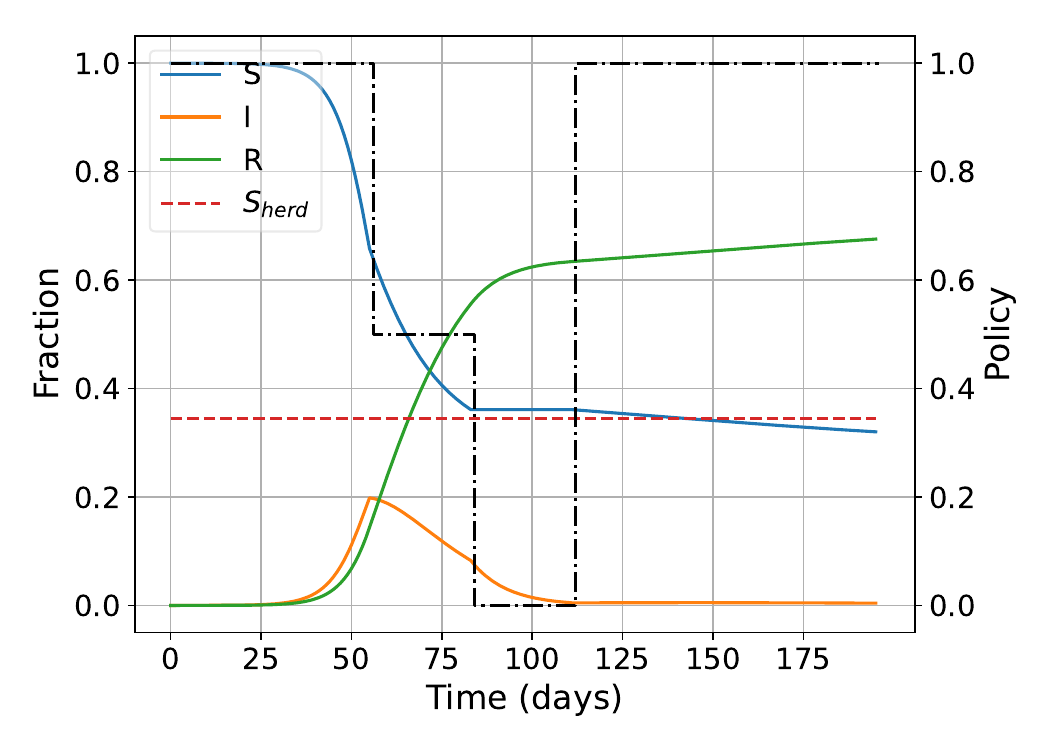}
         \caption{Optimal policy with the minimal policy time interval $\Delta t= 28$ days, $S_{\infty} = 0.32$}
         \label{fig:singcounty_dt=20}
    \end{subfigure}
    \hspace{1em}
    \begin{subfigure}[t]{0.45\linewidth}
    \centering
     \includegraphics[width=\linewidth]{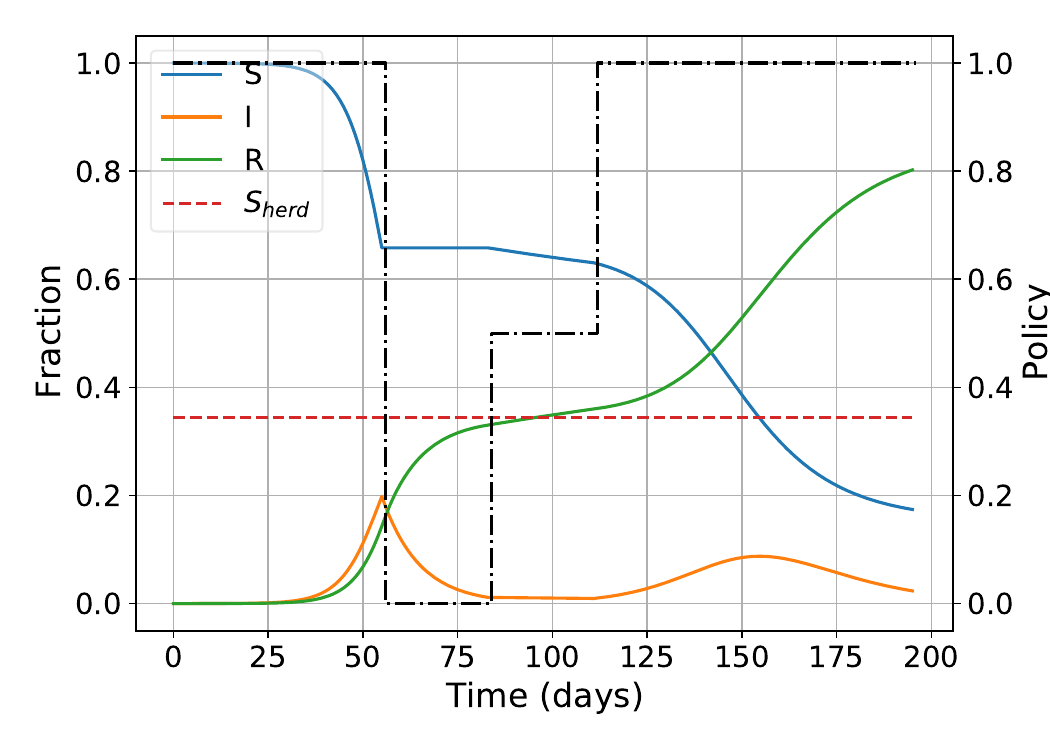}
     \caption{Flipped policy from panel (c), $S_{\infty} = 0.174$ }
     \label{fig:singcounty_dt=20_flipped}
    \end{subfigure}
\caption{Optimal policy and the SIR model of France from March 17 to May 11 2020. $S_{\text{herd}} \approx 0.345$.}
\label{fig:France}
\end{figure}

We compare the results from \cite{bliman_2021_optimal} to our model with the same cost function but only three possible levels of policy intensity $\alpha$. As in \cite{bliman_2021_optimal}, the general cost function \eqref{eqn:total_cost} reduces to the impact cost and is parameterized as the final epidemic size $R_{\infty}$. Bliman et al. assume that the paths considered all reach herd immunity. Therefore, in our search for the optimal policy, we exclude cases that do not reach herd immunity. Note that without this exclusion, the optimal solution is to adopt and hold the strictest possible policy starting from the beginning of the pandemic. This results in the least number of infections. For ease of computation, we consider three levels of policy intensity: 0, 0.5, 1 and fixed time intervals for the MPTI. We use the same set of parameters for the SIR model as in Bliman et al. \cite{bliman_2021_optimal}: $N = 6.7\times 10^{7}, I_0= 10^{3}$, $S_0= N-I_0$, $\mathcal{R}_0 = 2.9$.  Following \cite{bliman_2021_optimal}, we also choose 
the policy end time $T_0$ as close as possible to 100, thus setting $T_0 = 98$ since the time interval needs to be a multiple of the MPTI of 7 days. We show the result our algorithm produces in Fig.~\ref{fig:singcounty_dt=5} which we visually compare to the result from \cite{bliman_2021_optimal}, shown in Fig.~\ref{fig:bliman}. Note that we normalized curves by the total population. Both solutions are bang-bang controllers.  The solution using our model starts the control on day 63 (a multiple of 7) rather than day 61.9 (continuous). Slightly more people are infected under a  policy that is forced to use seven day intervals compared with continuous time as used by Bliman et al.  

Using a larger minimal policy time interval of 28 days and $T_0=112$, the optimal solution is no longer a bang-bang controller, as shown in Fig.~\ref{fig:singcounty_dt=20} with a larger $S_{\infty} = 0.32$. The optimal policy starts with a looser ``intermediate" policy phase followed by a stricter phase. Interestingly, in practice, during COVID-19 it was common for policies to start with the strictest restrictions followed by partial opening \cite{wikipedia_2022, department_of_public_health}. Thus, it is interesting to contrast the optimal policy with a policy in which the two stages are flipped in time, see Fig.~\ref{fig:singcounty_dt=20_flipped}. The flipped policy is a sub-optimal solution---it results in a larger pandemic size and a second wave of infections, as was often seen during the first two years of the COVID-19 pandemic.  Nevertheless, the policy in Fig.~\ref{fig:singcounty_dt=20_flipped}, while infecting more people, divides the impacted population into two distinctive waves, which could decrease daily hospital demand over the course of the outbreak.  Our policy model does not optimize for hospital demand.  Since many public health agencies (including Los Angeles County) considered hospital demand when making policy decisions, it could be important to consider in future studies.

\begin{table}[!h]
    \centering
    \begin{tabular}{||c|c|c|c|c|c|c|c||}
        \hline
       Figure & $T_0$ & $\Delta t$ &$N$ &$I_0$& $S_0$&$\mathcal{R}_0$&$S_{s\infty}$ \\
        \hline
        \hline
        A & 98 & 7 & $6.7\times 10^{7}$ &$10^{3}$ &$N-I_0$ &2.9 & 0.296\\
        \hline
        B & 100 & Not applicable & $6.7\times 10^{7}$ &$10^{3}$ &$N-I_0$ &2.9 & 0.31\\
        \hline
        C & 112 & 28 & $6.7\times 10^{7}$ &$10^{3}$ &$N-I_0$ &2.9 & 0.32\\
         \hline
        D & 112 & 28 & $6.7\times 10^{7}$ &$10^{3}$ &$N-I_0$ &2.9 & 0.174\\
        \hline
    \end{tabular}
    \caption{Parameters. }
    \label{tab:fig2para}
\end{table}

\subsubsection{Case Study - second wave in Los Angeles}
\begin{figure}[!ht]
    \begin{subfigure}[t]{0.48\textwidth}
        \centering
         \includegraphics[width=\textwidth]{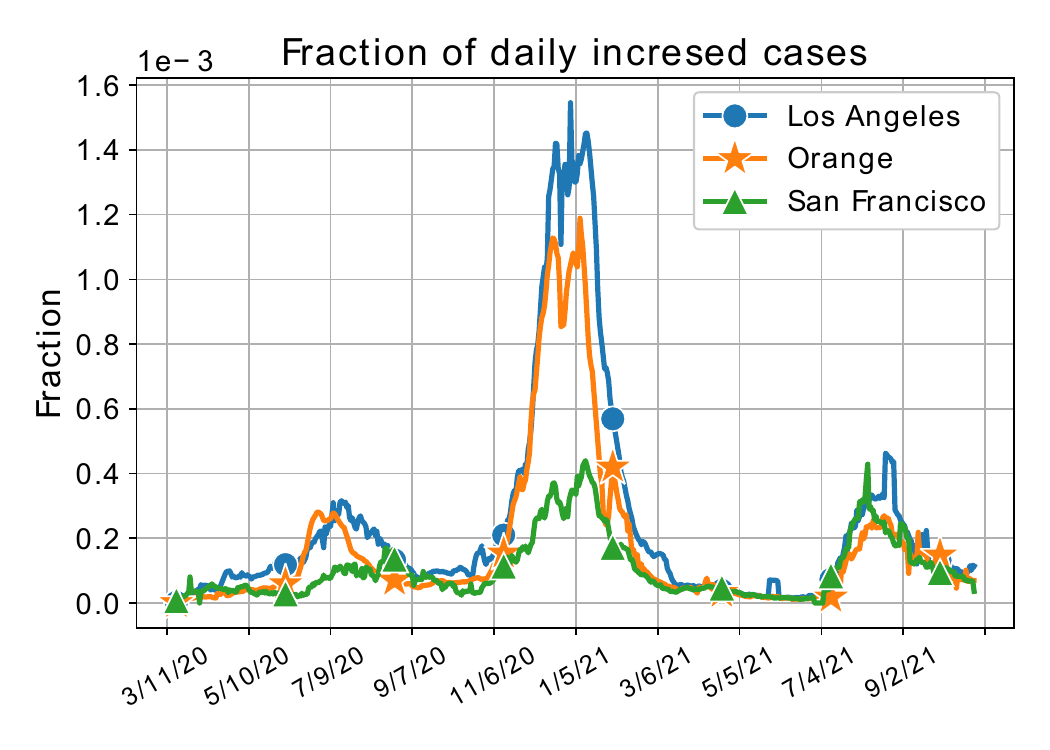}
         \caption{The fraction of the daily increase of the infected with a 7-day rolling average.}
         \label{fig:ca_i_real}
    \end{subfigure}
    \begin{subfigure}[t]{0.48\textwidth}
        \centering
         \includegraphics[width=\textwidth]{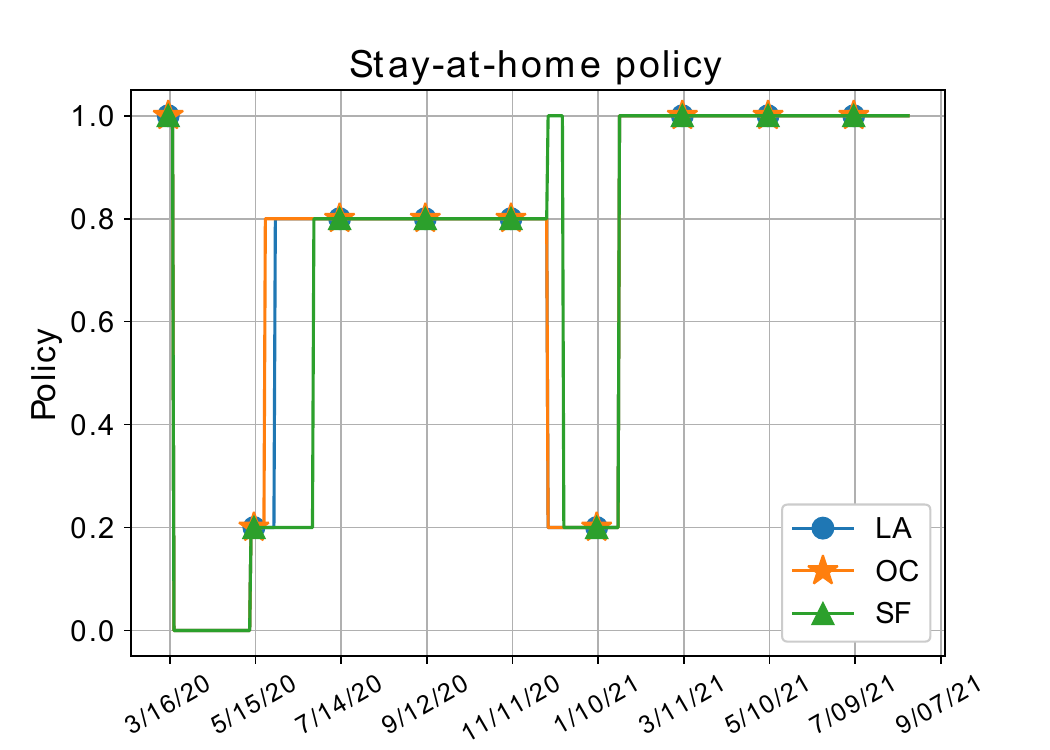}
         \caption{Stay-at-home Policy}
         \label{fig:ca_p_real}
    \end{subfigure}
    \caption{The fraction of the infected and `stay-at-home' policy over time in Los Angeles, San Francisco, and Orange County.}
    \label{fig:ca_real}
\end{figure}

We first present the course of infections in three counties in California and their corresponding `stay-at-home' policy changes from Mar 2020 to Sept 2021.
Fig. \ref{fig:ca_i_real} shows the 7-day rolling average of the fraction of the daily increased infected cases based on the data from \cite{dong2020interactive} in 3 counties with the largest population density in California, namely, San Francisco, Orange, and Los Angeles. There were 3 major outbreaks during the given time interval. For the first and second waves, Orange and Los Angeles Counties followed similar trajectories, while San Francisco County stayed more contained. Due to substantial holiday travel in winter 2020-21, the second wave was a much larger than the first.\\
In \cite{us_policy}, the US Centers for Disease Control and Prevention describes six levels of `stay-at-home' policy.  The intensity of the policy decreases as the numerical value increases. The exact descriptions of the five levels of policies and their numerical representation are shown in Table \ref{tab:policy_table}. Fig.~\ref{fig:ca_p_real} shows the change of the intensity of the `stay-at-home' policy during the same period.  Policy during the first wave was  \textit{proactive}, whereas for the second wave it was more \textit{reactive}. 
This may reflect some hesitancy on the part of policy-makers as well as lesser compliance by the population at large by the time the second wave emerged. During the second wave, with a relatively strict policy, the regions all stayed below herd immunity. With vaccination available in early 2021, the pandemic in all three regions tapered off. 
\begin{table}[!h]
    \centering
    \scalebox{0.8}{
    \begin{tabular}{||c|c||}
        \hline
        Numerical value & `Stay-at-home' policy\\
        \hline
        \hline
        0 & Mandatory for all individuals \\
        \hline
        0.2 & Mandatory only for all individuals in certain areas of the jurisdiction \\
        \hline
        0.4 & Mandatory only for at-risk individuals in the jurisdiction \\
        \hline
        0.6 & Mandatory only for at-risk individuals in certain areas of the jurisdiction\\
        \hline
        0.8 & Advisory/Recommendation\\
         \hline
        1 & No order for individuals to stay home\\
         \hline
    \end{tabular}
    }
    \caption{CDC stay-at-home policies. There are 6 levels of policies and we map the levels linearly onto the interval $[0,1]$ for simplicity. The numerical value on the left is used to graph actual policies over time in Fig.~\ref{fig:ca_p_real}.
    }
    \label{tab:policy_table}
    \vspace{-4mm}
\end{table}

Now we consider a counterfactual study of how the pandemic would have evolved had herd immunity been reached during the second wave, controlled by our policy model, using parameters measured from the Los Angeles data. We choose to study the period of the second wave for several reasons. First, the data reporting scheme improved for the second wave compared to the first wave. In addition, with the experience and knowledge gained from the first wave, authorities were in a better position to make optimal decisions. Given that there was no complete lockdown during the second wave, we consider the policy intensity levels $A = \{0.2, 0.6, 1\}$, and use the minimal policy time interval $\Delta t = 7$. We choose $0.2$ as our maximal policy intensity because full lockdown was not desirable during this period.  We choose a second policy level of $0.6$ as the midpoint between $0.2$ and $1$.  In all simulations we optimize for final pandemic size and compare the optimal controls found.

\begin{figure}[!th]
    \centering
    \begin{subfigure}[b]{0.48\textwidth}
         \includegraphics[width=\textwidth]{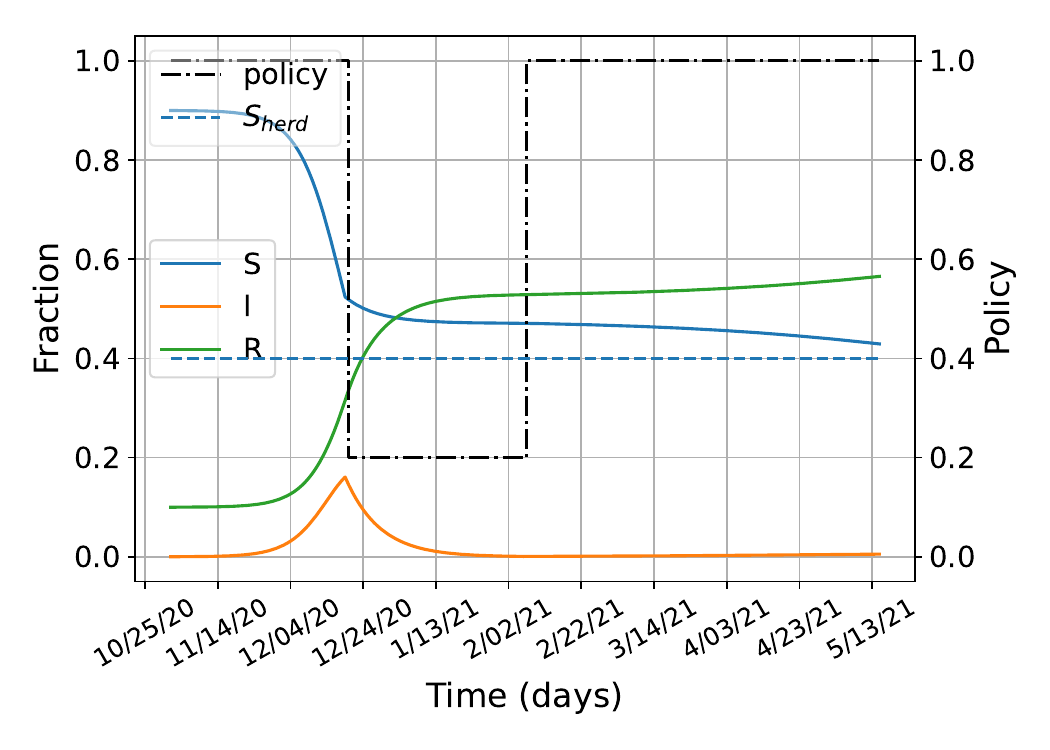}
         \caption{$R_0=2.5$, $r_0 = 0.1$.}
         \label{fig:la_2.5_0.1}
    \end{subfigure}
    \begin{subfigure}[b]{0.48\textwidth}
     \includegraphics[width=\textwidth]{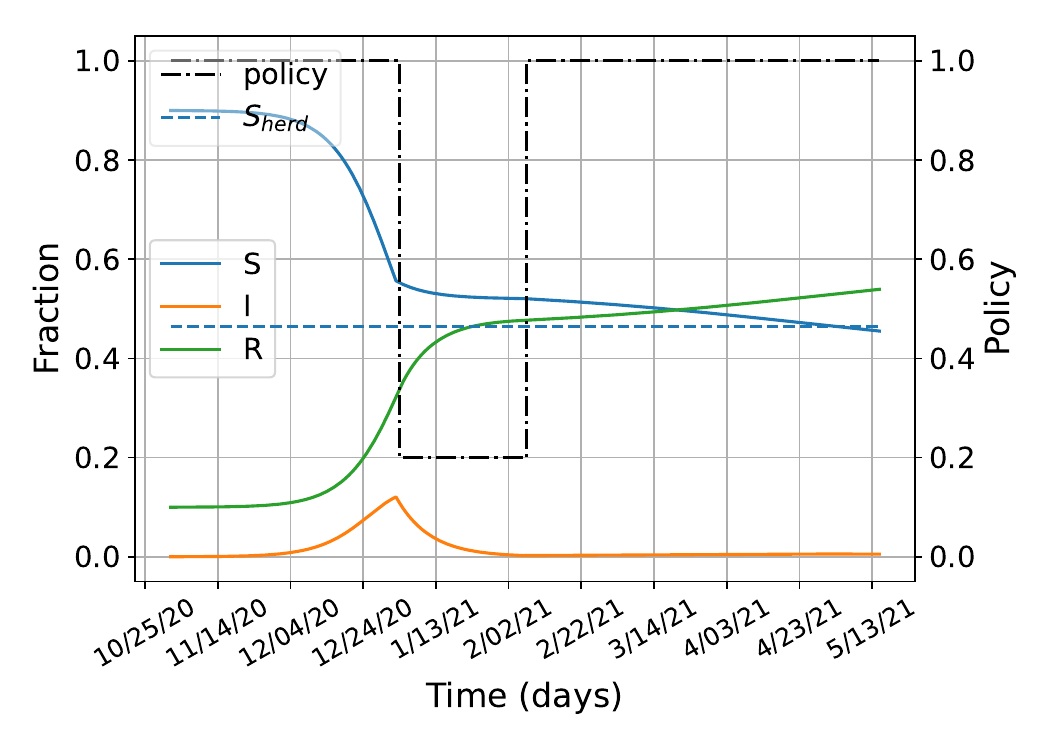}
     \caption{$R_0=2.15$, $r_0 = 0.1$.}
     \label{fig:la_2.15_0.1}
    \end{subfigure}
    \\
    \begin{subfigure}[b]{0.48\textwidth}
        \centering
         \includegraphics[width=\textwidth]{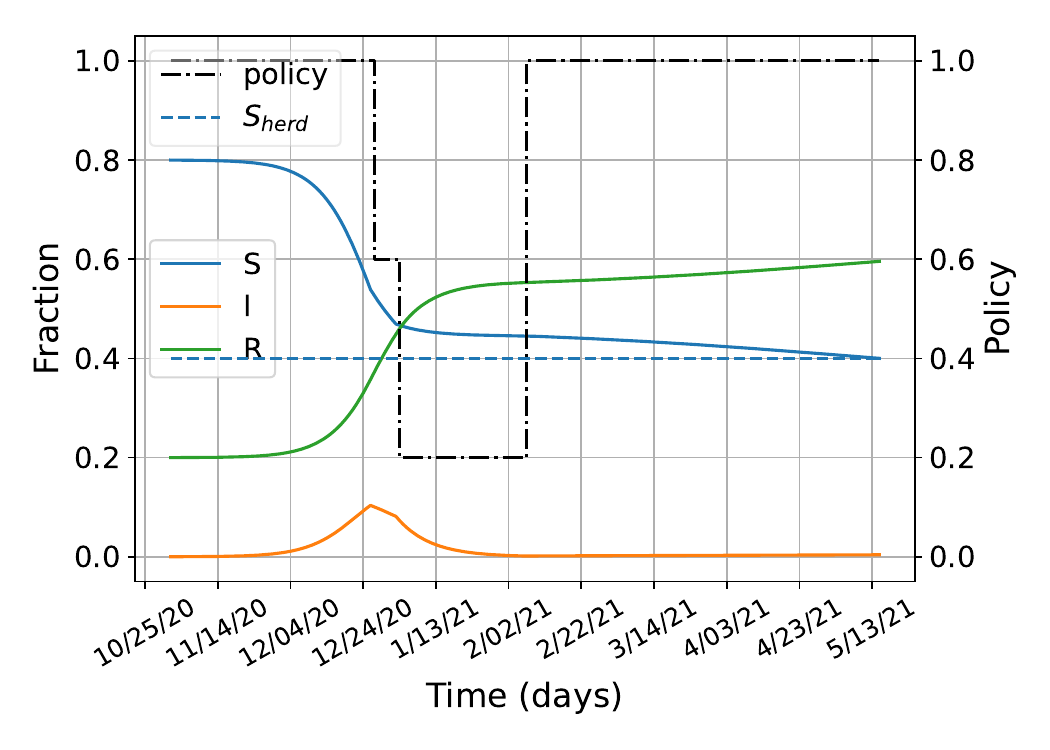}
         \caption{$R_0=2.5$, $r_0 = 0.2$}
         \label{fig:la_2.5_0.2}
    \end{subfigure}
    \begin{subfigure}[b]{0.48\textwidth}
        \centering
         \includegraphics[width=\textwidth]{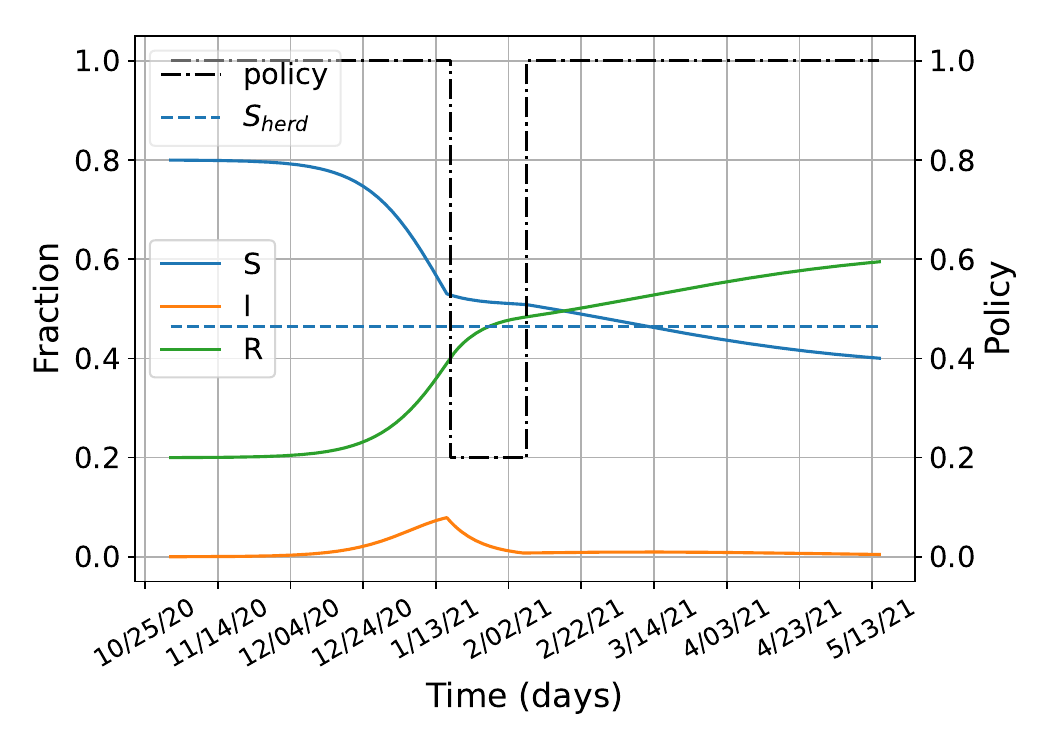}
         \caption{$R_0=2.15$, $r_0 = 0.2$}
         \label{fig:la_2.15_0.2}
    \end{subfigure}
    \\
        \begin{subfigure}[b]{0.48\textwidth}
        \centering
         \includegraphics[width=\textwidth]{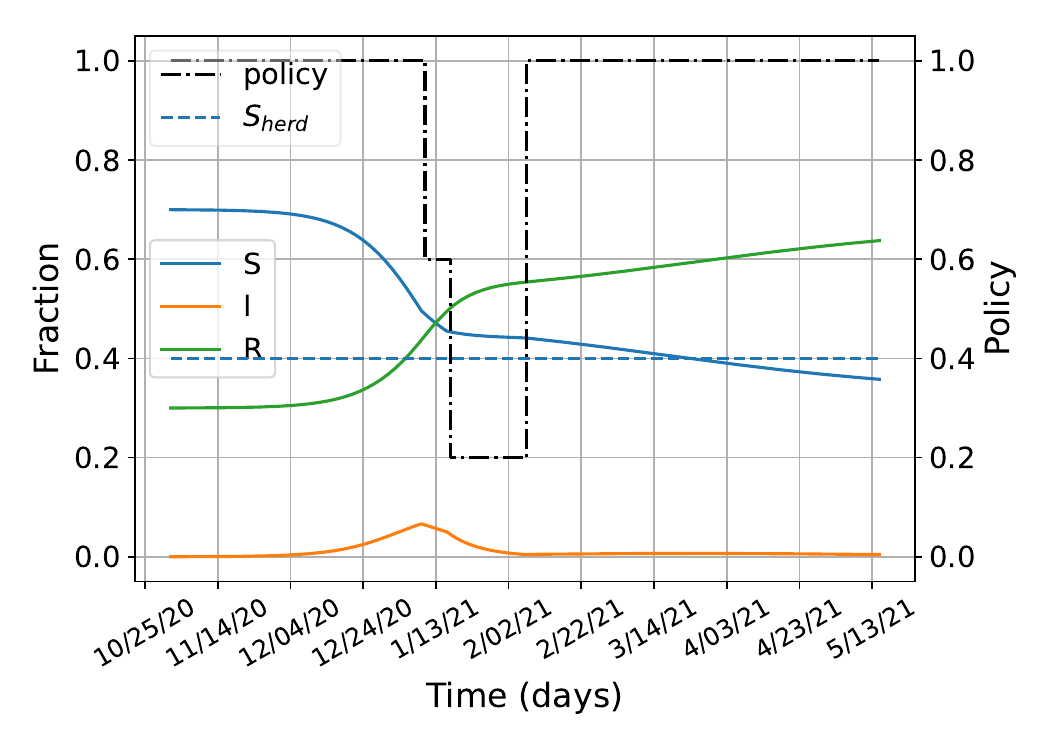}
         \caption{$R_0=2.5$, $r_0 = 0.3$}
         \label{fig:la_2.5_0.3}
    \end{subfigure}
    \begin{subfigure}[b]{0.48\textwidth}
        \centering
         \includegraphics[width=\textwidth]{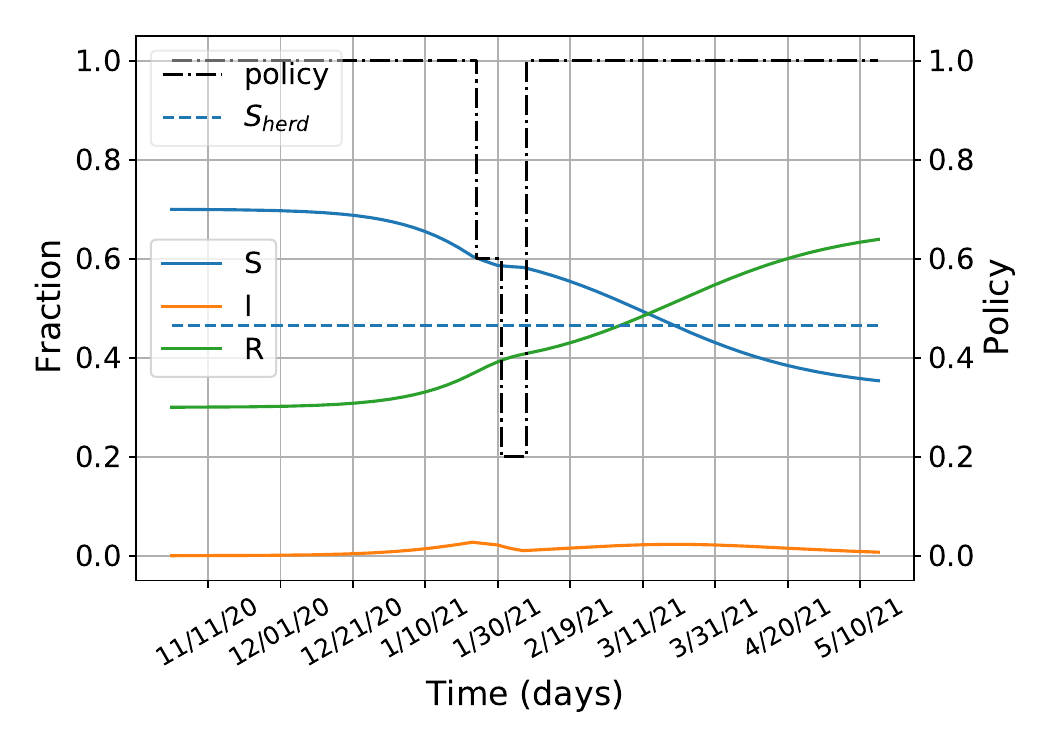}
         \caption{$R_0=2.15$, $r_0 = 0.3$.}
         \label{fig:la_2.15_0.3}
    \end{subfigure}
\caption{Optimal policy in Los Angeles with the basic reproduction number $R_0 = 2.5, 2.15$, $S_{\text{herd}} \approx 0.4, 0.465$,  and the fraction of the initial recovered population $r_0 = 0.1,0.2,0.3$, respectively}
\label{fig:la_sim}
\end{figure}

In Fig. \ref{fig:la_sim}, the left column (Figs.~\ref{fig:la_2.5_0.1}, \ref{fig:la_2.5_0.2}, \ref{fig:la_2.5_0.3}) is the simulated SIR with the optimal policy when the basic reproduction number $R_0 = 2.5$ and the initial recovered $r_0 = 0.1,0.2,0.3$. The right column (Figs.~\ref{fig:la_2.15_0.1}, \ref{fig:la_2.15_0.2}, \ref{fig:la_2.15_0.3}) is the simulated SIR with the optimal policy when the reproduction number $R_0 = 2.15$ and the initial recovered $r_0 = 0.1,0.2,0.3$. This value $R_0=2.5$ is estimated from the early COVID-19 infected data (Jan 22—Mar 15, 2020 \cite{dong2020interactive}) and $R_0=2.15$ is estimated using the infected data from September 16 to November 15, 2021 (\cite{dong2020interactive}), prior to the second wave. All optimal policies have a bang-bang-like shape. The policy started approximately around the peak of the infected curve, and the resulting dynamics approach herd immunity. For larger values of $r_0$, we expect that a shorter period of high intensity policy is needed to reach herd immunity and our results confirm this. Once enough of the population is infected and recovered, a shorter control policy is needed.

\section{Multi-layer multi-regional case}
\begin{figure}[!h]
    \centering
    \includegraphics[width=0.65\linewidth]{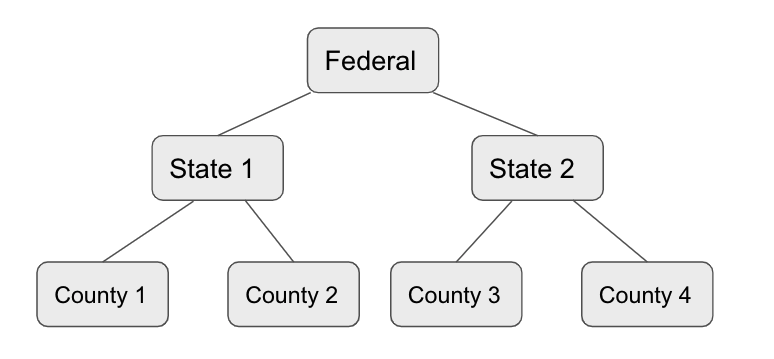}
    \caption{An example of a three-layer hierarchical structure.}
    \label{fig:hierachical}
\end{figure}

\begin{figure}
    \centering
    \includegraphics[width=0.6\linewidth]{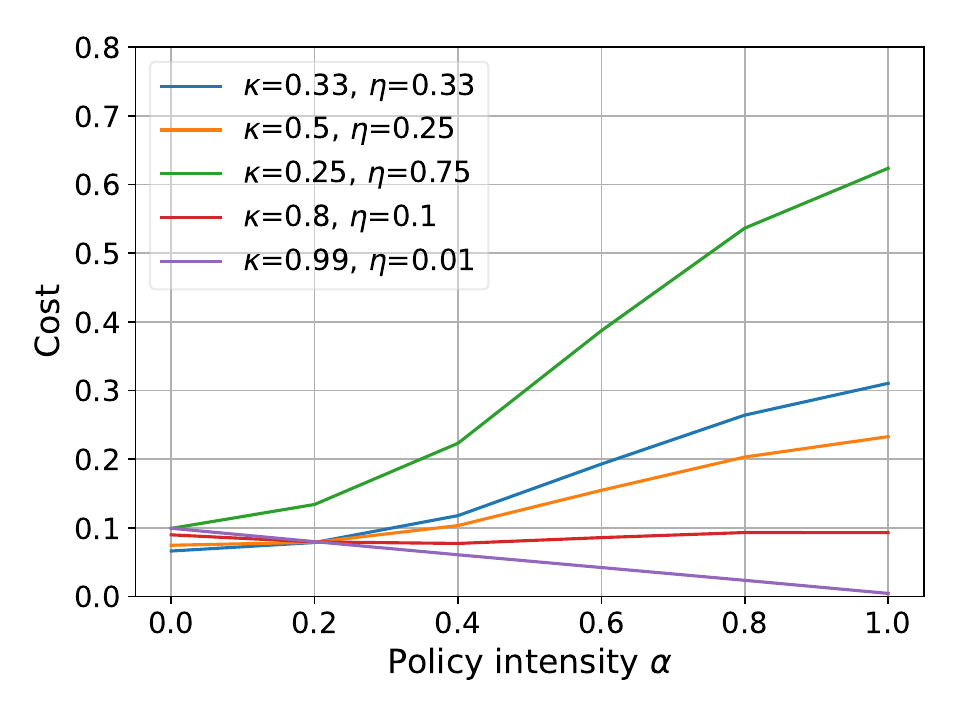}
    \caption{Different cost functions vs policy intensity $\alpha$.}
    \label{fig:cost}
\end{figure}

In this section, we present a multi-regional model with multiple policy-making layers, extending the model propose by Jia {\em et al} \cite{jia2021game}  to consider a dynamic epidemic model and control policies discussed above (see Section \ref{sec:single_region}).
Specifically, we propose a game-theoretic model in which regions are combined into layers, with the top layer corresponding to the highest-level decision maker (e.g., a federal government), the next layer comprised of the next level of decision making (e.g., states or provinces), and so on (see Fig. \ref{fig:hierachical}).
The top decision maker chooses the policy first, next all the decision makers in the next layer simultaneously, and so on.
Additionally, we consider a special case in which there are multiple decision makers (e.g., states, counties, etc) choosing their epidemic control policies simultaneously in one layer. We use a form of hierarchical best response dynamic to compute approximate equilibria in this multi-layer game~\cite{jia2021game}, performing this computation independently for each time interval (essentially assuming that the players do not reason explicitly about future dynamics when making instantaneous policy decisions at a particular point in time).

The multi-region case naturally has a competition between regions to optimize their strategy with respect to the choices made by other regions.  For this reason the single-region model does not directly extend.  
There are two main differences between our work and that of Jia at al. \cite{jia2021game}. First, their model is based on Nash equilibrium, where agents make decisions with other agents' possible actions in mind. We use the idea of `in-game learning' \cite{learning-in-games}. 
We assume that the agents gradually evolve towards the best decisions instead of being optimal instantly.
In practice, each region in the game assumes other regions' policies (at the same level) stay the same when optimizing its own cost function. Second, we focus on the dynamics, instead of a snapshot in time considered by Jia et al. 
\paragraph{Network SIR}
In practice, counties can hardly be treated as independent. People travel across county borders to work and socialize. The majority of the literature of network-style SIR models focus on the individuals as nodes and study the effects of interpersonal network on the pandemics \cite{network_book,1556-1801_2022_3_359, netowrk_spreaders}. For example, \cite{netowrk_spreaders} empirically study how well various centrality measures perform at identifying which individuals in a network will be the best spreaders of disease. In \cite{pnas_linear}, the authors explains why most COVID-19 infection curves are linear after the first peak in the context of the contact network using a network SIR model. There are a few works that study the interplay between different geographical regions rather than the interpersonal contact network. In \cite{pnas_kernel_sir}, a kernel-modulated SIR model was introduced to model the spread of COVID-19 across counties. The kernel is based on the spatial proximity between regions. Metapopulation epidemic models are based on the spatial structure of the environment, and detailed knowledge of transportation infrastructure and movement patterns. The metapopulation dynamics of infectious diseases has generated a wealth of models and results using mechanistic approaches taking explicitly into account the movement of individuals (\cite{grais2003assessing,keeling2002estimating,sars2006}). For example, in \cite{sars2006}, the authors proposed a multi-regional compartmental model using medical geography theory (central place theory) and studied the effect of the travel of individuals (especially those infected and exposed) between regions on the global spread of severe acute respiratory syndrome (SARS).
Another way to account for the interplay between regions is to use a cross excitation matrix \cite{hawkes_baichuan}. This scheme assumes the a uniform mixing of the population across regions and the infected population in one region can trigger the infection in another. The entries of the matrix records the pair-wise cross excitation from one region to another. 
In this paper, we assume uniform mixing in the population and use an excitation matrix $K = \{K_{aa'}\}$ to model the travel and infections across counties. Our network-style SIR is the following:
\begin{eqnarray}\displaystyle
\begin{cases}
&\dfrac{d S_a(t)}{dt} =  - \alpha_a \beta \sum_{a'}K_{aa'}\frac{I_{a'}(t)S_a(t)}{N_a} , \\
&\dfrac{dI_a(t)}{dt} = \alpha_a\beta \sum_{a'}K_{aa'}\frac{I_{a'}(t)S(t)}{N_a} -\gamma  I_a(t) , \\
&\dfrac{dR_a(t)}{dt} =  \gamma  I_a(t) , \\
&  S(0)=    S_0,  \quad  I(0)=  I_0, \quad   R(0)=  R_0.
\end{cases}
\end{eqnarray}
For any county $a$, the rate of change from $S_a$ to $I_a$ triggered by $I_{a'}$ depends on $K_{aa'}$,   the current fraction of the susceptible $S_a$ in county $a$ and the current fraction of the infected $I_{a'}$ in county $a'$. Note that $K_{aa} = 1$. When $K = I$, the network SIR is the independent SIR.
\paragraph{Cost function}
Consider the $i$-th time interval $[i\Delta t,(i+1)\Delta t]$ and $u(t) = \alpha$ for $t\in[i\Delta t,(i+1)\Delta t]$. A Region $a$ adopts the following cost function:
\[
c_{i\Delta t,(i+1)\Delta t}^a(\alpha) = \kappa_{a}(1-\alpha)\Delta t/T + \eta_{a} R_{a,T}(\alpha) + (1-\kappa_{a}-\eta_{a})(\alpha - \pi (\alpha))^2\Delta t/T,
\]
where $R_{a,T}(\alpha)$ is the epidemic size of region $a$ at time $T$. For a top-layer region $f$, there is no non-compliance cost and the cost function is
\[
c_{\Delta t}^f(\alpha) = \kappa_f(1-\alpha)\Delta t/T + \eta_f R_{f,T}(\alpha),
\]
where $\kappa_f + \eta_f = 1$ and $R_{f,T}(\alpha)$ is the number of the recovered of region $f$ which is an aggregation of the the epidemic size of its leaf nodes. 


\subsection{Algorithms}
The single region algorithm minimizes over all admissible piece-wise functions, while the multiple-region algorithm only minimizes over every time interval. We assume there are up to three layers: federal government, the states, and the counties. 
At $n$-th time interval, we first determine the optimal policy intensity that minimizes the cost $C^f_{n\Delta t,(n+1)\Delta t}$ for the federal layer. After obtaining the optimal federal policy, each state optimizes its own cost function $C^s_{n\Delta t,(n+1)\Delta t}$ for the period $[n\Delta t, (n+1)\Delta t]$ unilaterally, i.e., assuming other states follow their previous policies. Next, we choose the optimal policy intensity for the counties in the same manner. Note that the federal layer does not pay the non-compliance cost as it is not subject to any higher-level policy making. The states and counties may pay a non-compliance cost. The full details of the three-layer model is in Alg.~\ref{alg:gamesir}.
\begin{algorithm}[!h]
\caption{\textsc{Game Policy SIR}}\label{alg:gamesir}
\begin{algorithmic}[1]
    \STATE \textbf{Input}: \textbf{Time $T$},\textbf{ excitation matrix} $K$, \textbf{ intensity levels} $A$,  \textbf{time interval $\Delta t$}
    \STATE Initialize \textbf{state, county policies}
    \STATE Number of policy stages $N = \frac{T}{\Delta t}$, $n = 1$
    \WHILE{$n \leq N$}
        \STATE $t = n\Delta t$
         \WHILE{$t < T$}
             \FOR{every state $s$}
                \FOR{every county $a$ in state $s$}
                    \STATE update $S_a,I_a,R_a$ according to the current policy $\alpha_a$ and the excitation matrix $K$:
                    \STATE $S_a(t) =  S_a(t-1)- \alpha_a \beta \sum_{a'}K_{aa'}\frac{I_{a'}(t-1)S_a(t-1)}{N_a}$\\
                    \STATE $I_a(t) = I_a(t-1) +  \alpha_a\beta \sum_{a'}K_{aa'}\frac{I_{a'}(t-1)S(t-1)}{N_a} -\gamma  I_a(t-1)$\\
                    \STATE $R_a(t) = R_a(t-1) + \gamma  I_a(t-1)$
                \ENDFOR
            \ENDFOR
            \STATE $t \mathrel{+}= 1$
        \ENDWHILE
        \STATE $\alpha_f = \arg\min_{\alpha' \in A} C^f_{n\Delta t,(n+1)\Delta t}(\alpha')$
        \FOR{every state $s$}
            \STATE $\alpha_s = \arg\min_{\alpha' \in A} C^s_{n\Delta t,(n+1)\Delta t}(\alpha')$
            \FOR{every county $a$ in state $s$}
            \STATE $\alpha_a = \arg\min_{\alpha' \in A} C^a_{n\Delta t,(n+1)\Delta t}(\alpha')$
            \ENDFOR
        \ENDFOR
         \STATE $n \mathrel{+}= 1$
    \ENDWHILE
\end{algorithmic}
\end{algorithm}

\subsection{Simulations}
In this section, we present results for a three-county example of the multiple regions game and a three-county example with a state. First, we discuss when one layer exists (i.e., only counties). The game between the counties is through cross excitation of infection among the counties. Next, we study the case when a governing state is added.

We consider three interacting counties with the excitation matrix $K$:
$$
K = \begin{bmatrix}
1 & 0 & 0\\
0.1 & 1 & 0\\
0 & 0.1 & 1
\end{bmatrix}
$$
We set the reproduction number $\mathcal{R}_0 = 2$ and therefore, $S_{\text{herd}} = 0.5$. Counties 1, 2, 3 have initial fractions of the infected population as $i_0 = 0.2, 0.1,0.1$, respectively.  This implies that county 1 has a bigger outbreak initially, and part of the infection in county 2 is excited from county 1 and part of the infection in county 3 is excited from county 2. The cost functions for all counties consist of an implementation cost and an impact cost with equal weights ($\eta_a=\kappa_a = 1/2$, for all $a$). The minimal policy time interval $\Delta$ is set to be $7$ (days).

\begin{figure}[!ht]
    \centering
    \begin{subfigure}[t]{0.45\linewidth}
        \centering
         \includegraphics[width=\linewidth]{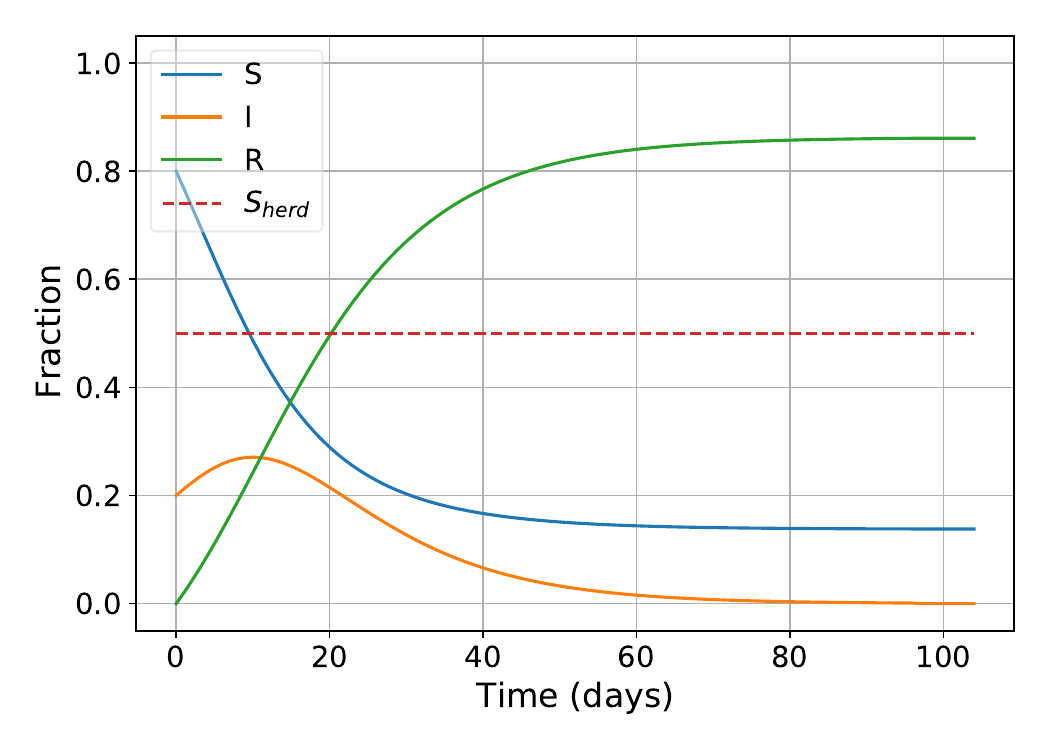}
         \caption{County 1. No intervention.}
         \label{fig:c1_w/o}
    \end{subfigure}
    \begin{subfigure}[t]{0.45\linewidth}
    \centering
     \includegraphics[width=\linewidth]{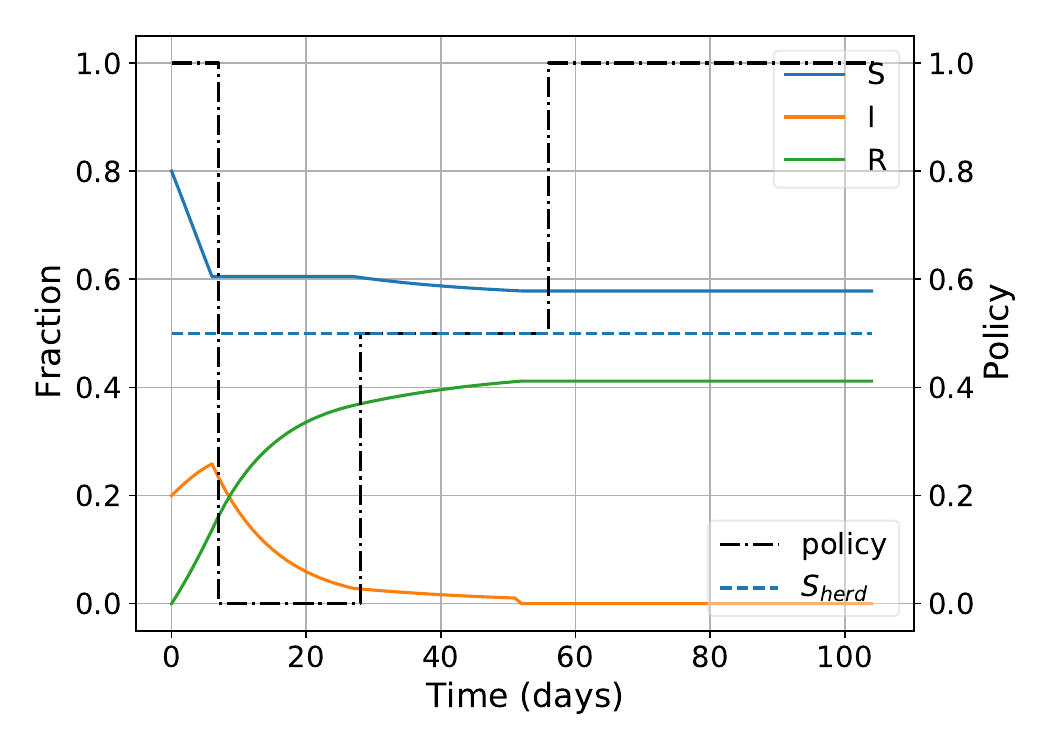}
     \caption{County 1.}
     \label{fig:c1_w}
    \end{subfigure}
    \\
    \begin{subfigure}[b]{0.45\textwidth}
        \centering
         \includegraphics[width=\textwidth]{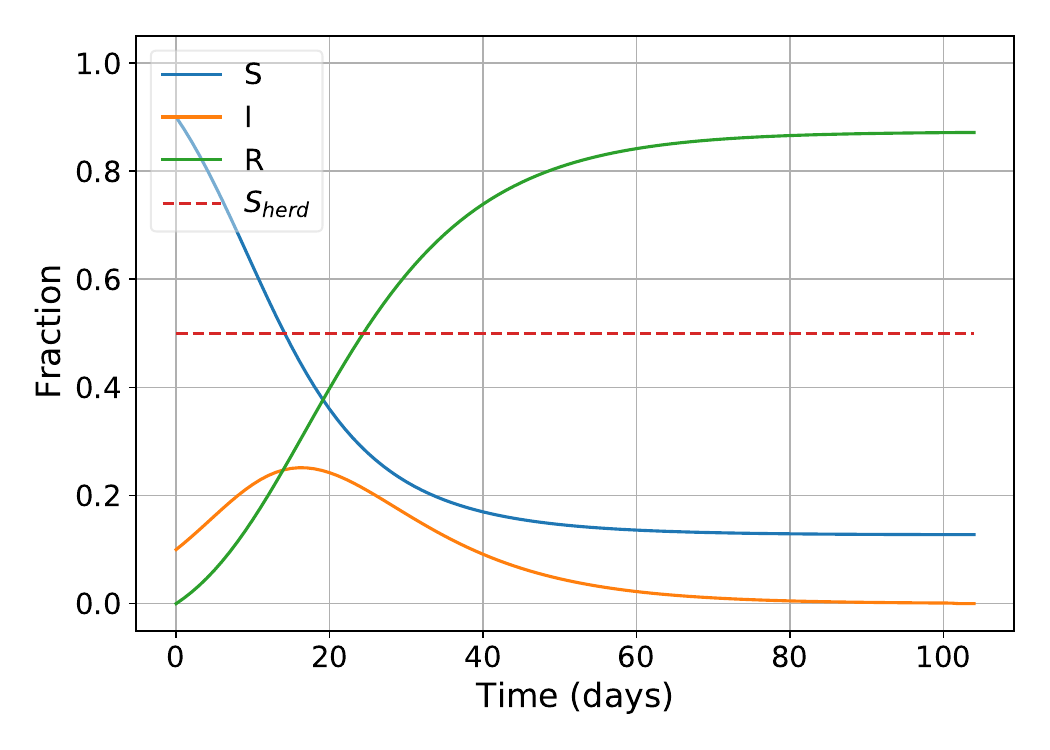}
         \caption{County 2. No intervention.}
         \label{fig:c2_w/o}
    \end{subfigure}
    \begin{subfigure}[b]{0.45\textwidth}
        \centering
         \includegraphics[width=\textwidth]{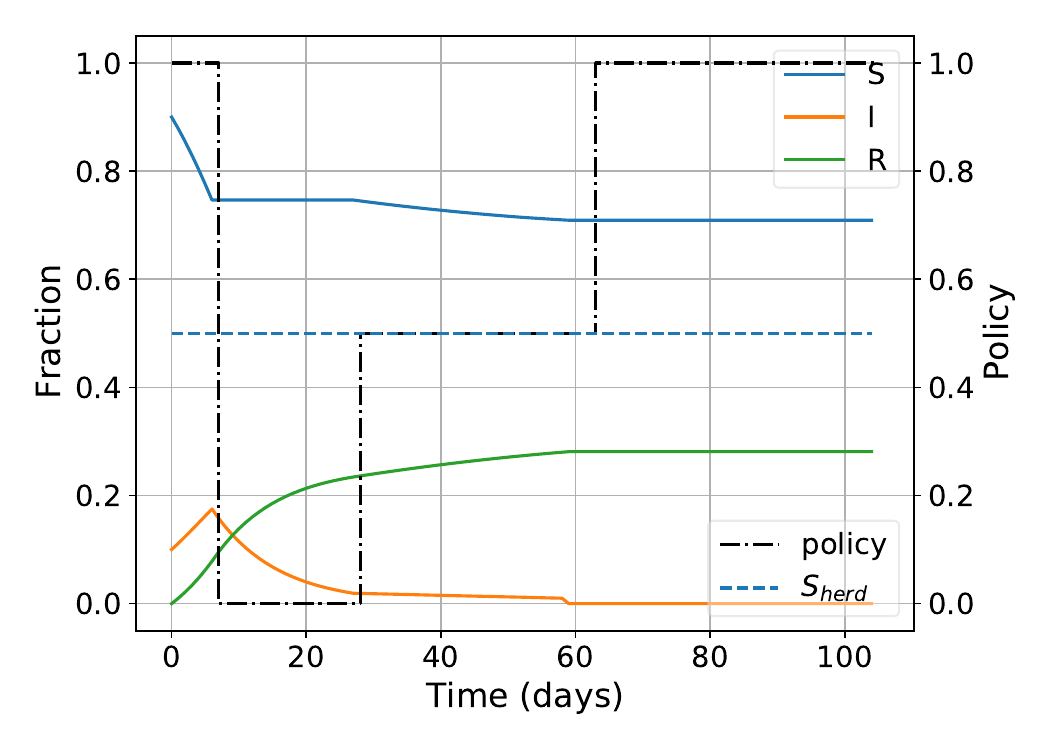}
         \caption{County 2.}
         \label{fig:c2_w}
    \end{subfigure}
    \\
    \begin{subfigure}[b]{0.45\textwidth}
        \centering
         \includegraphics[width=\textwidth]{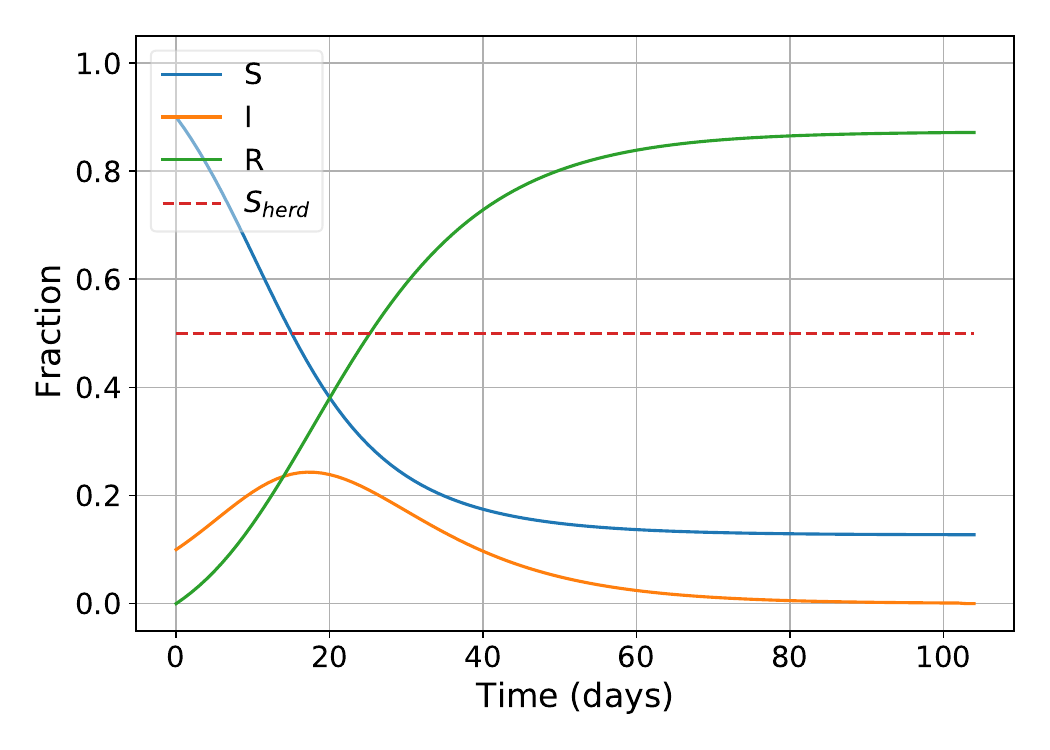}
         \caption{County 3. No intervention.}
         \label{fig:c3_w/o}
    \end{subfigure}
    \begin{subfigure}[b]{0.45\textwidth}
        \centering
         \includegraphics[width=\textwidth]{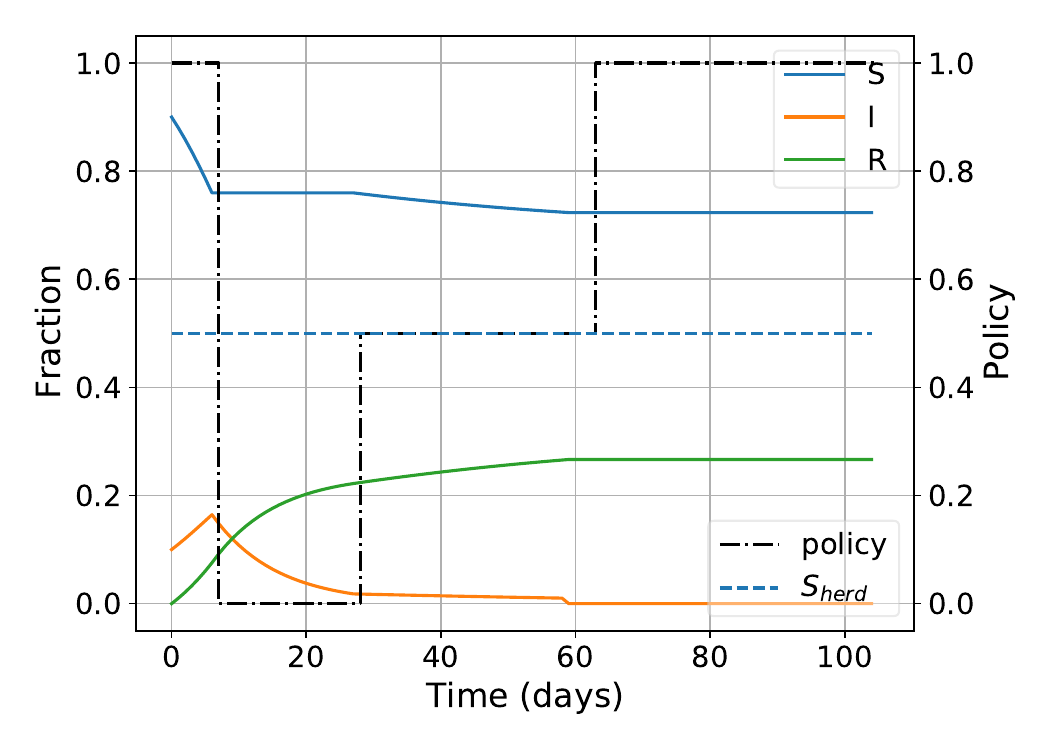}
         \caption{County 3.}
         \label{fig:c3_w}
    \end{subfigure}
\caption{An example of three dependent counties without and with interventions. With intervention, for all counties, the coefficients for the implementation cost $\kappa = \frac{1}{2}$ and the coefficients for the impact cost $\eta = \frac{1}{2}$. The minimal policy time interval $\Delta t = 7$.}
\label{fig:3-regions}
\end{figure}

\begin{figure}[!ht]
    \centering
    \begin{subfigure}[b]{0.45\linewidth}
        \centering
         \includegraphics[width=\linewidth]{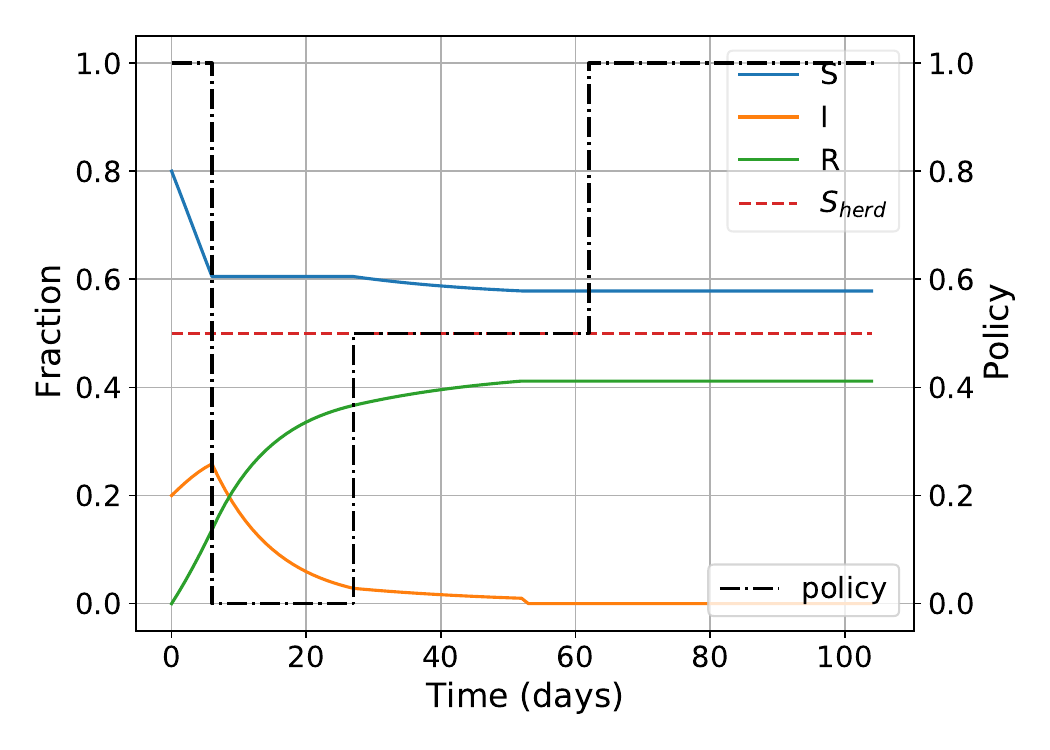}
         \caption{County 1 with a state.}
         \label{fig:sc1}
    \end{subfigure}
    \begin{subfigure}[b]{0.45\textwidth}
        \centering
         \includegraphics[width=\textwidth]{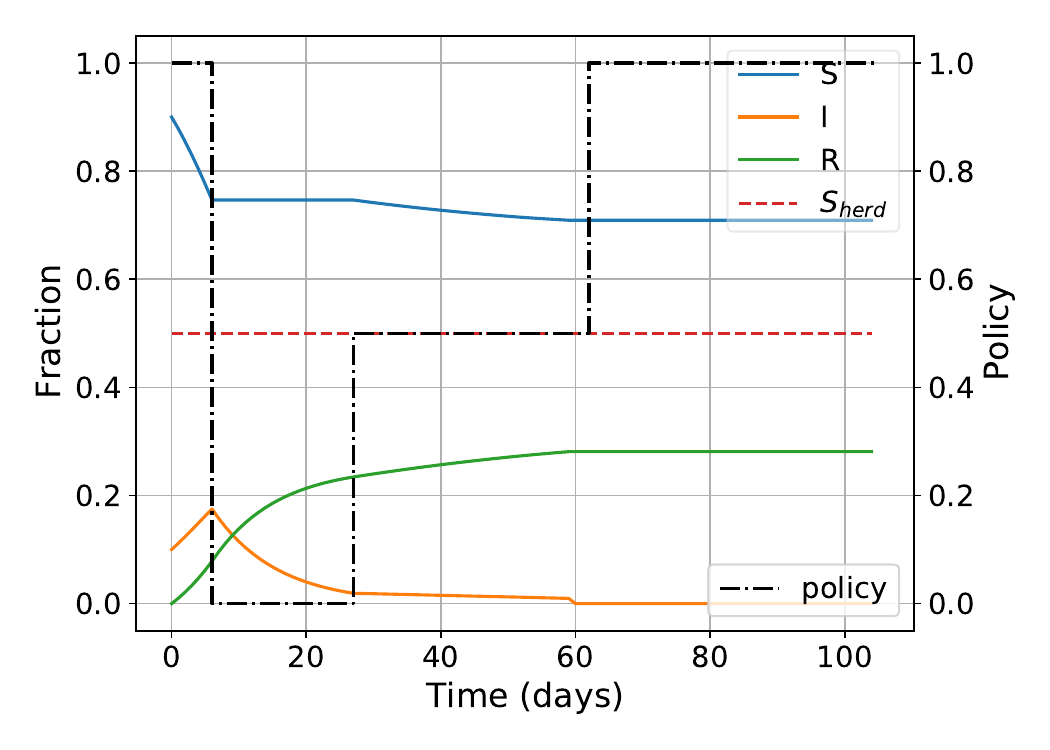}
         \caption{County 2 with a state.}
         \label{fig:sc2}
    \end{subfigure}
    \\
    \begin{subfigure}[b]{0.45\textwidth}
        \centering
         \includegraphics[width=\textwidth]{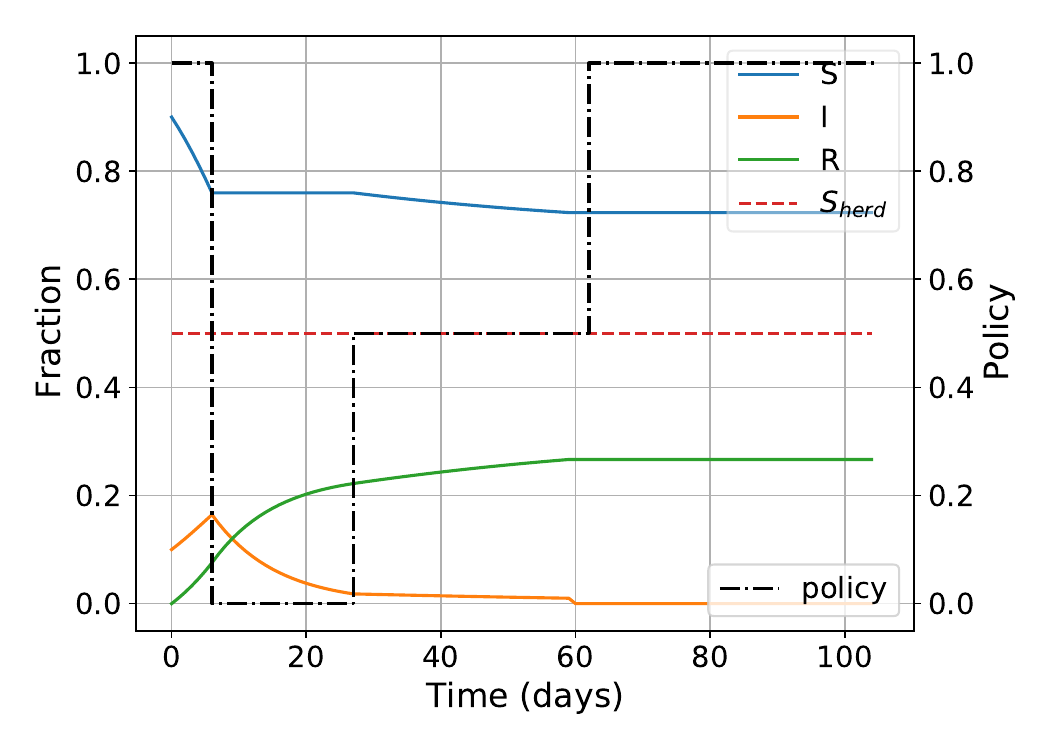}
         \caption{County 3 with a state.}
         \label{fig:sc3}
     \end{subfigure}
    \begin{subfigure}[b]{0.45\linewidth}
        \centering
         \includegraphics[width=\linewidth]{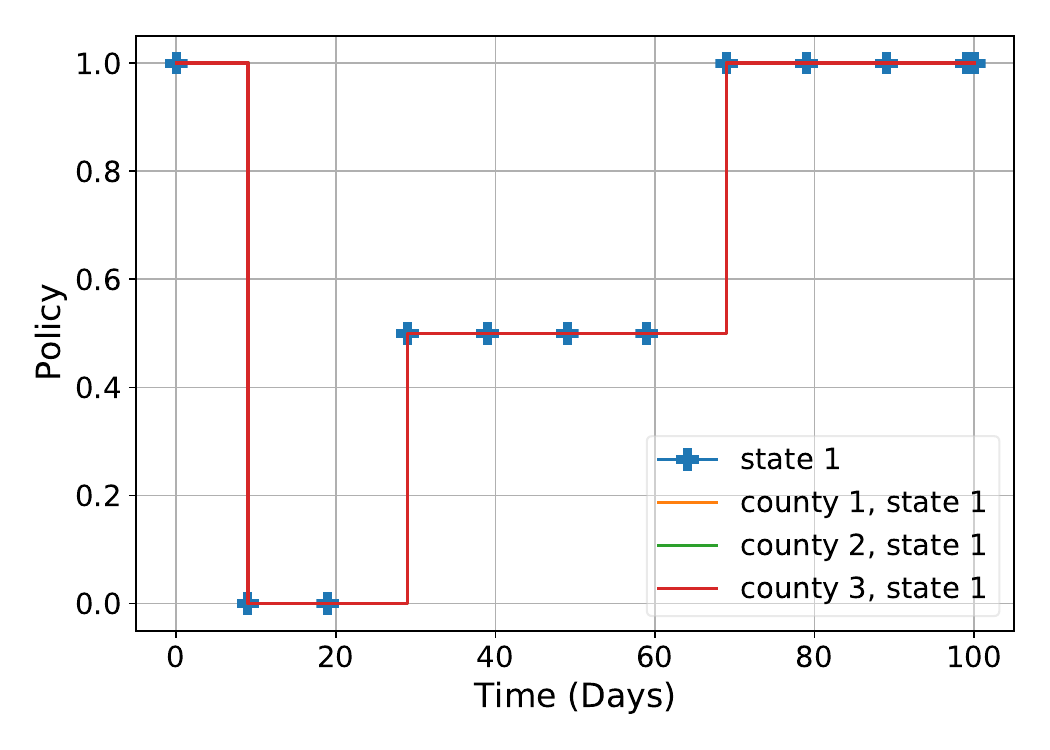}
         \caption{State policy and counties' policies.}
         \label{fig:sc_policy}
        \end{subfigure}
    \caption{An example with 3 counties and a governing state. For all counties, the coefficients for the implementation cost $\kappa = \frac{1}{6}$, the coefficients for the impact cost $\eta = \frac{1}{6}$ and the coefficients for the non-compliance cost $1-\kappa-\eta = \frac{2}{3}$. For the state, the coefficients for the implementation cost $\kappa = \frac{1}{3}$, the coefficients for the impact cost $\eta = \frac{2}{3}$. The minimal policy time interval $\Delta t = 7$.}
\label{fig:sc}
\end{figure}

The left column (Figs. \ref{fig:c1_w/o}, \ref{fig:c2_w/o}, \ref{fig:c3_w/o}) are simulations for the counties without any intervention and the right column  (Figs. \ref{fig:c1_w}, \ref{fig:c2_w}, \ref{fig:c3_w}) are simulations with interventions. Without intervention, we see propagation of waves of infection from county 1 to county 2 and then to county 3. All of the counties reached herd immunity eventually. With interventions, policy restrictions started on day $7$ and, for county 2 and 3, the infected curves decrease before reaching their peaks. With control, county 1 contained the pandemic and the final $S_{\infty}$ is close herd immunity level $S_{herd}$. With a fewer infected population to begin with, county 2 and 3 contained the pandemic before reaching herd immunity.
Fig. \ref{fig:sc} shows the results of adding a governing state on top of the county layer. We keep the ratio of the weights for the implementation cost and the impact cost to be 1:1, the same as in the no-state case in Fig. \ref{fig:3-regions}. The state has slightly different weights, with the ratio of the weights for the implementation cost and the impact cost being 1:2. Compared to Fig. \ref{fig:3-regions}, by adding a state, the three counties ended up with the same policy. In this case, the noncompliance cost results in each county choosing the same policy as the state rather than different policies.

\section{Discussion and future work}
We propose a policy-making model coupled with the SIR model to study a single region and game-like interactions between multiple regions. The model demonstrates its ability to model real-life situations with different sets of parameters in both one and multiple-region scenarios. One can extend the model to a hierarchical structure by building multiple layers of the multiple regions model and study the cross-layer effects.

In the search for an optimal policy, we used a naive depth-first search algorithm for the one-region model. One can speed up the algorithm by removing some of the obvious non-optimal paths.

In our model, the policy intensity $\alpha$ is a heuristic representation of the lockdown, social distancing and mask policy. It remains to be discussed how other policies, for example, vaccination policies, affects the spreading in the different stages of the pandemic. The model ignores some of the important features like the limitation of the hospital capacity \cite{thompson2020will}, which could be added as constraints when minimizing the cost function. Fig.~\ref{fig:ca_p_real} shows the policy for the first wave is proactive while the one for the second wave is reactive. One possible effect is from fatigue of following policy, which increases in time and has a memory. So far, the model does not have the capability of modeling this fatigue. In the future, one could consider an adaptive term in the cost function to model it.  The network example considered was rather simplistic, with just three counties within one state.  One could consider more complex systems with multiple layers.  The computational method here would likely need to be improved to address the computational complexity of the search space. In addition, a potentially important generalization is to capture implementation and impact costs with more refined cross-layer dependencies, but this is potentially non-trivial from both a modeling and computational perspective.


\newpage
\vfill\eject
\bibliographystyle{plain}
\bibliography{refer.bib}

\begin{thebibliography}{10}

\bibitem{office_city_of_los_angeles}
{COVID}-19: Keeping {L}os {A}ngeles safe.
\newblock \url{https://coronavirus.lacity.org/english/mayor-garcetti-issues-safer-home-emergency-order-stopping-non-essential-activities-outside}, 2020.
\newblock Accessed: 2022-09-06.

\bibitem{sd_reduce_fz}
Pierre-Alexandre Bliman and Michel Duprez.
\newblock How best can finite-time social distancing reduce epidemic final size?
\newblock {\em Journal of Theoretical Biology}, 511:110557, 2021.

\bibitem{bliman_2021_optimal}
Pierre-Alexandre Bliman, Michel Duprez, Yannick Privat, and Nicolas Vauchelet.
\newblock Optimal immunity control and final size minimization by social distancing for the {SIR} epidemic model.
\newblock {\em Journal of Optimization Theory and Applications}, 03 2021.

\bibitem{blower2002health}
SM~Blower, Katia Koelle, and John Mills.
\newblock Health policy modeling: Epidemic control, {HIV} vaccines, and risky behavior.
\newblock {\em Quantitative evaluation of HIV prevention programs}, pages 260--289, 2002.

\bibitem{pnas_1918}
Martin C.~J. Bootsma and Neil~M. Ferguson.
\newblock The effect of public health measures on the 1918 influenza pandemic in {U.S.} cities.
\newblock {\em Proceedings of the National Academy of Sciences}, 104(18):7588--7593, 2007.

\bibitem{burke2006}
Donald~S. Burke, Joshua~M. Epstein, Derek~A.T. Cummings, Jon~I. Parker, Kenneth~C. Cline, Ramesh~M. Singa, and Shubha Chakravarty.
\newblock Individual-based computational modeling of smallpox epidemic control strategies.
\newblock {\em Academic Emergency Medicine}, 13(11):1142--1149, 2006.

\bibitem{CDC_timeline}
{Centers for Disease Control and Prevention}.
\newblock {CDC Museum {COVID}-19 Timeline}.
\newblock \href{https://www.cdc.gov/museum/timeline/covid19.html}{https://www.cdc.gov/museum/timeline/covid19.html}.
\newblock Accessed: 2022-06-20.

\bibitem{us_policy}
{Centers for Disease Control and Prevention}.
\newblock {U.S. State and Territorial Stay-At-Home Orders: March 15, 2020 – August 15, 2021 by County by Day}.
\newblock \href{https://data.cdc.gov/Policy-Surveillance/U-S-State-and-Territorial-Stay-At-Home-Orders-Marc/y2iy-8irm}{https://data.cdc.gov/Policy-Surveillance/U-S-State-and-Territorial-Stay-At-Home-Orders-Marc/y2iy-8irm}.
\newblock Accessed: 2022-06-20.

\bibitem{department_of_public_health}
{Department of Public Health, Los Angeles}.
\newblock {Department of Public Health, Los Angeles}.
\newblock \href{http://publichealth.lacounty.gov/}{http://publichealth.lacounty.gov/}, 3 2022.
\newblock Accessed: 2022-06-20.

\bibitem{dong2020interactive}
Ensheng Dong, Hongru Du, and Lauren Gardner.
\newblock An interactive web-based dashboard to track {{COVID}}-19 in real time.
\newblock {\em The Lancet Infectious Diseases}, 2020.

\bibitem{learning-in-games}
Drew Fudenberg and David~K. Levine.
\newblock {\em {The Theory of Learning in Games}}, volume~1 of {\em MIT Press Books}.
\newblock The MIT Press, 1998.

\bibitem{pnas_kernel_sir}
Xiaolong Geng, Gabriel~G. Katul, Firas Gerges, Elie Bou-Zeid, Hani Nassif, and Michel~C. Boufadel.
\newblock A kernel-modulated {SIR} model for {COVID}-19 contagious spread from county to continent.
\newblock {\em Proceedings of the National Academy of Sciences}, 118(21):e2023321118, 2021.

\bibitem{grais2003assessing}
Rebecca~F Grais, J~Hugh Ellis, and Gregory~E Glass.
\newblock Assessing the impact of airline travel on the geographic spread of pandemic influenza.
\newblock {\em European journal of epidemiology}, 18(11):1065--1072, 2003.

\bibitem{jia2021game}
Feiran Jia, Aditya Mate, Zun Li, Shahin Jabbari, Mithun Chakraborty, Milind Tambe, Michael Wellman, and Yevgeniy Vorobeychik.
\newblock A game-theoretic approach for hierarchical policy-making.
\newblock {\em arXiv preprint arXiv:2102.10646}, 2021.

\bibitem{keeling2002estimating}
Matt~J Keeling and Pejman Rohani.
\newblock Estimating spatial coupling in epidemiological systems: a mechanistic approach.
\newblock {\em Ecology Letters}, 5(1):20--29, 2002.

\bibitem{network_book}
Istvan~Z Kiss, Joel~C Miller, and Peter Simon.
\newblock {\em (Book) Mathematics of Epidemics on Networks: from Exact to Approximate Models}.
\newblock Springer, 2017.

\bibitem{1556-1801_2022_3_359}
Prateek Kunwar, Oleksandr Markovichenko, Monique Chyba, Yuriy Mileyko, Alice Koniges, and Thomas Lee.
\newblock A study of computational and conceptual complexities of compartment and agent based models.
\newblock {\em Networks and Heterogeneous Media}, 17(3):359--384, 2022.

\bibitem{ijerph19106119}
Thomas~H. Lee, Bobby Do, Levi Dantzinger, Joshua Holmes, Monique Chyba, Steven Hankins, Edward Mersereau, Kenneth Hara, and Victoria~Y. Fan.
\newblock Mitigation planning and policies informed by {COVID}-19 modeling: A framework and case study of the state of {H}awaii.
\newblock {\em International Journal of Environmental Research and Public Health}, 19(10), 2022.

\bibitem{lin_2020}
Rong-Gong Lin~II.
\newblock L.a. county now requires residents to wear face coverings. here are the details, Apr 2020.

\bibitem{netowrk_spreaders}
Brian Macdonald, Paulo Shakarian, Nicholas Howard, and Geoffrey Moores.
\newblock Spreaders in the network {SIR} model: an empirical study.
\newblock {\em arXiv preprint arXiv:1208.4269}, 2012.

\bibitem{Allison2018}
Allison~L. Pitt, Keith Humphreys, and Margaret~L. Brandeau.
\newblock Modeling health benefits and harms of public policy responses to the us opioid epidemic.
\newblock {\em American Journal of Public Health}, 108(10):1394--1400, 2018.
\newblock PMID: 30138057.

\bibitem{sars2006}
Shigui Ruan, Wendi Wang, and Simon~A. Levin.
\newblock The effect of global travel on the spread of {SARS}.
\newblock {\em Mathematical Biosciences and Engineering}, 3(1):205--218, 2006.

\bibitem{dfs}
Robert Tarjan.
\newblock Depth-first search and linear graph algorithms.
\newblock In {\em 12th Annual Symposium on Switching and Automata Theory (swat 1971)}, pages 114--121, 1971.

\bibitem{thompson2020will}
RN~Thompson, CA~Gilligan, and NJ~Cunniffe.
\newblock Will an outbreak exceed available resources for control? {E}stimating the risk from invading pathogens using practical definitions of a severe epidemic.
\newblock {\em Journal of the Royal Society Interface}, 17(172):20200690, 2020.

\bibitem{pnas_linear}
Stefan Thurner, Peter Klimek, and Rudolf Hanel.
\newblock A network-based explanation of why most {COVID}-19 infection curves are linear.
\newblock {\em Proceedings of the National Academy of Sciences}, 117(37):22684--22689, 2020.

\bibitem{ca_timeline}
Wikipedia.
\newblock {California {COVID}-19 Timeline}.
\newblock \href{https://en.wikipedia.org/wiki/Timeline\_of\_the\_COVID-19\_pandemic\_in\_California}{https://en.wikipedia.org/wiki/Timeline\_of\_the\_COVID-19\_pandemic\_in\_California}.
\newblock Accessed: 2022-06-20.

\bibitem{wikipedia_2022}
Wikipedia.
\newblock {COVID}-19 pandemic in {F}rance.
\newblock \url{https://en.wikipedia.org/wiki/COVID-19_pandemic_in_France#Timeline_of_measures}.
\newblock Accessed: 2022-09-06.

\bibitem{hawkes_baichuan}
Baichuan Yuan, Hao Li, Andrea~L. Bertozzi, P.~Jeffrey Brantingham, and Mason~A. Porter.
\newblock Multivariate spatiotemporal {H}awkes processes and network reconstruction.
\newblock {\em SIAM Journal on Mathematics of Data Science}, 1(2):356--382, 2019.

\end{thebibliography}

\end{document}